\documentclass[12pt]{article}
\usepackage{amssymb, amsmath, amsthm}
\usepackage[all]{xy}

\renewcommand{\baselinestretch}{1.1}

\newcommand{\mr}[1]{\mathrm{#1}}
\newcommand{\mf}[1]{\mathfrak{#1}}
\newcommand{\mc}[1]{\mathcal{#1}}

\newcommand{\Z}{{\bf Z}}
\newcommand{\Q}{{\bf Q}}
\newcommand{\zp}{{\bf Z}_p}
\newcommand{\qp}{{\bf Q}_p}

\newcommand{\II}{\mf{I}}
\newcommand{\I}{\mc{I}}

\newcommand{\tlam}{\tilde{\Lambda}}
\newcommand{\tgam}{\tilde{\Gamma}}
\newcommand{\tdel}{\tilde{\Delta}}

\newcommand{\eq}[1]{(\ref{#1})}

\DeclareMathOperator{\Hom}{Hom} \DeclareMathOperator{\Aut}{Aut}
\DeclareMathOperator{\End}{End} \DeclareMathOperator{\Gal}{Gal}
\DeclareMathOperator{\Sel}{Sel} 
\DeclareMathOperator{\coker}{coker}

\newtheorem{theorem}{Theorem}[section]
\newtheorem{proposition}[theorem]{Proposition}
\newtheorem{lemma}[theorem]{Lemma}
\newtheorem{corollary}[theorem]{Corollary}

\newtheorem{hypothesis}{Hypothesis}

\newtheorem{conjecture}[theorem]{Conjecture}

\theoremstyle{definition}

\theoremstyle{remark}
\newtheorem*{remark}{Remark}
\newtheorem*{example}{Example}
\newtheorem*{ack}{Acknowledgments}

\numberwithin{equation}{section}

\begin{document}

\title{Iwasawa theory and the Eisenstein ideal}
\author{Romyar T. Sharifi}
\date{}
\maketitle
\renewcommand{\baselinestretch}{1}
\begin{abstract}
	We verify, for each odd prime $p < 1000$, a conjecture of W.\ McCallum and
	the author's on the surjectivity of pairings on $p$-units constructed out 
	of the cup product on the first Galois cohomology group of the maximal
	unramified outside $p$ extension of $\Q(\mu_p)$ with $\mu_p$-coefficients.  
	In the course of the proof, we relate several Iwasawa-theoretic and 
	Hida-theoretic objects.  In particular, we construct 
	a canonical isomorphism between an Eisenstein ideal modulo its square and
	the second graded piece in an augmentation filtration of a classical Iwasawa
	module over an abelian pro-$p$ Kummer extension of the cyclotomic 
	$\zp$-extension of an abelian field. This Kummer extension arises from the 
	Galois representation on an inverse limit of ordinary parts of first
	cohomology groups of modular curves that was considered by M.\ Ohta in order
	to give another proof of the Iwasawa Main Conjecture in the spirit of 
	that of
	B.\ Mazur and A.\ Wiles.  In turn, we relate the Iwasawa module over the
	Kummer extension to the quotient of the tensor product of the classical
	cyclotomic Iwasawa module and the Galois group of the Kummer extension by
	the image of a certain reciprocity map that is constructed out of an 
	inverse limit of cup products up the cyclotomic tower.  We give an
	application to the structure of the Selmer groups of Ohta's modular 	
	representation taken modulo the Eisenstein ideal.
\end{abstract}
\renewcommand{\baselinestretch}{1.2}

\section{Introduction}

The central object of our study is, for an odd prime $p$, the cup product in the cohomology of the Galois group $\mc{G}$ of the maximal pro-$p$ unramified outside $p$ extension of $F = \Q(\mu_p)$ with $\mu_p$-coefficients,
$$
	H^1(\mc{G},\mu_p) \otimes H^1(\mc{G},\mu_p) \xrightarrow{\cup}
	H^2(\mc{G},\mu_p^{\otimes 2}).
$$
This cup product is well known to describe the structure of $\mc{G}$ modulo the closed normal subgroup generated by $p$th powers and triple commutators.  Moreover, it provides a useful first object in the study of Iwasawa theory over metabelian pro-$p$ unramified outside $p$ extensions of $F$, as detailed in \cite{me-massey} and Section \ref{cupproducts}. 

Via Kummer theory, the cup product induces a pairing on the $p$-units $\mc{E}_F$ of $F$,
$$
(\, \cdot \, , \, \cdot \,) \colon \mc{E}_F \times \mc{E}_F \to A_F \otimes \mu_p,
$$
where $A_F$ denotes the $p$-part of the class group of $F$.   In \cite{mcs}, W.\ McCallum and the author studied this pairing in detail.   For $p$ satisfying Vandiver's Conjecture, which states that $A_F$ equals its $(-1)$-eigenspace $A_F^-$ under complex conjugation,
we conjectured that the image of the pairing $(\, \cdot\, , \, \cdot \,)$ generates $A_F \otimes 
\mu_p$ \cite[Conjecture 5.3]{mcs}.   Denoting the projection of $(\,\cdot\, , \,\cdot\,)$ to $A_F^- \otimes \mu_p$ by $(\,\cdot\, , \,\cdot\,)^-$, we state the obvious unconditional strengthening of that conjecture.

\begin{conjecture}[McCallum-Sharifi] \label{pairingconj}
	The image of the pairing $(\,\cdot\, , \,\cdot\,)^-$ generates $A_F^- 
	\otimes \mu_p$.
\end{conjecture}

In fact, the author expects that $(\,\cdot\, , \,\cdot\,)^-$ is itself always surjective.
As described in \cite{mcs}, a positive solution of this conjecture has numerous applications: among others, to R.\ Greenberg's pseudo-nullity conjecture in Iwasawa theory, to products in $K$-groups of cyclotomic integer rings, and to Y.\ Ihara's pro-$p$ Lie algebra arising from the outer representation of Galois on the pro-$p$ fundamental group of the projective line minus three points.  In this paper, we focus instead on a relationship between the structure of certain Iwasawa modules over Kummer extensions and the structure of ordinary Hecke algebras of modular forms localized at the Eisenstein ideal.  It is this intriguing relationship that enables us to prove the following result.

\begin{theorem} \label{mainresult}
	The pairing $(\, \cdot \, , \,\cdot\,)$ is surjective for $p < 1000$. 
\end{theorem}

The cup product pairing has the property that $(x,1-x) = 0$ if $x$ and $1-x$ are both
$p$-units in $F$.
By imposing relations arising from this fact, McCallum and the author were able to compute a unique nontrivial possibility, up to scalar, for each projection of $(\,\cdot\, , \,\cdot\,)$ to a nontrivial $\Delta = \Gal(F/\Q)$-eigenspace of $A_F$ for all $p <$ 10,000 (now completed for $p <$ 25,000).  Therefore, we have a complete computation of the pairing $(\, \cdot\, , \,\cdot\,)$ up to a unit in $(\Z/p\Z)[\Delta]$ for all $p < 1000$.

Our method of proof, while given in the context of Iwasawa and Hida theory, is in the spirit of K.\ Ribet's proof of the converse of Herbrand's theorem \cite{ribet}.  We now describe the basic idea at the finite level, first recalling Ribet's work and then explaining our extension of it.  Let $\omega$ denote the Teichm\"uller character.  If $p$ divides the Bernoulli number $B_k$ for some even $k$ less than $p$, then there is a congruence between a weight $2$, level $p$, and character $\omega^{k-2}$ Eisenstein series and a cuspidal eigenform $f$ of the same weight, level, and character at a prime $\mf{p}$ above $p$ in the ring of coefficients of $f$.  Ribet showed that the fixed field of the kernel of the Galois representation $\rho_f$ attached to $f$ on $G_F$ contains a nontrivial unramified elementary abelian $p$-extension $H$ of $F$ such that $\Gal(H/F)$ has an $\omega^{1-k}$-action of $\Delta$.  In fact, 
$H$ is fixed by the kernel of the action of $G_F$ on a choice
of lattice for $\rho_f$ modulo $\mf{p}$.

For our construction, we assume that $p$ does not divide the Bernoulli number $B_{p+1-k}$ and, for simplicity of this description, that $p^2$ does not divide $B_k$.  We find a metabelian extension $M$ of $F$ in the fixed field of the kernel of $\rho_f$, the arithmetic of which determines the nontriviality of a certain value of the cup product.  
Unlike $H$, the field $M$ cannot be found in the fixed field of the representation on a lattice modulo $\mf{p}$.  To find it, we pass to a larger quotient of a lattice
by a sublattice intermediate between those determined by $\mf{p}$ and $\mf{p}^2$
(see Section \ref{galrep}).  
As an extension of $F$, the field $M$ has maximal abelian subextension equal to
the compositum $HE$ of $H$ and a degree $p$ unramified outside $p$ extension $E = F(\eta^{1/p})$ (with $\eta \in \mc{E}_F$) of $F$ that has an $\omega^{k-1}$-action of $\Delta$.  Over $HE$, the field $M$ is an unramified elementary abelian $p$-extension. 
We show that if no prime above $p$ splits completely in $M/HE$, then $(p,\eta)$ does not vanish.  To show that no prime above $p$ splits completely in $M/HE$ is 
equivalent to showing that the $p$th coefficient of $f$ is not $1$ modulo $\mf{p}^2$ (for some $f$ as above), which is computable.

We will consider a rather general set of abelian fields and Dirichlet characters in this article, but let us continue to describe the situation over $K = \Q(\mu_{p^{\infty}})$ here.
Let $X_K$ denote the Galois group of the maximal unramified abelian pro-$p$ extension of $K$.  
The Iwasawa Main Conjecture states that the $\omega^{1-k}$-eigenspace of $X_K$ has characteristic ideal generated by a $p$-adic power series determined by $L_p(\omega^k,s)$.  This conjecture was first proven by B.\ Mazur and A.\ Wiles in \cite{mw} via a study of the action of the absolute Galois group on ``good" quotients of the Jacobians of certain modular curves.  The basic idea of their proof is the same as that for Ribet's theorem, but the proof is vastly more involved.   A similar but slightly more streamlined proof was given more recently by M.\ Ohta \cite{ohta2} using Hida theory (see also \cite{wiles}).  Ohta's proof extends methods of M.\ Kurihara \cite{kurihara} and G.\ Harder and R.\ Pink \cite{hp}, who developed the analogue of Ribet's work in the more technical setting of level one modular forms.  It is Ohta's work that we shall call upon in this paper.

We study the action of the absolute Galois group $G_{\Q}$ on an Eisenstein component $\mf{X}$ of a particular inverse limit of cohomology groups of modular curves 
previously considered by Ohta \cite{ohta-comp}
(see Section \ref{prelim} for the precise definition).  Define $\mf{h}$ to be the localization of Hida's ordinary cuspidal Hecke algebra at the maximal ideal containing the Eisenstein ideal $\I$ with character $\omega^k$.  Then $\mf{X}$ is an $\mf{h}$-module.  If $p$ divides $B_k$ but not $B_{p+1-k}$, Ohta showed that $\mf{X}$ is free of rank $2$ over $\mf{h}$.  
The fixed field of the kernel of the representation $\rho \colon G_{\Q} \to \mr{Aut}_{\mf{h}}\, \mf{X}$ contains a certain abelian pro-$p$ extension $L$ of $K$ that is unramified outside $p$ and totally ramified above $p$.  We can consider $X_L$, the Galois group of the maximal unramified abelian pro-$p$ extension of $L$, and its filtration by powers of the augmentation ideal $I_G$ in $\zp[[G]]$, where $G = \Gal(L/K)$.  

Using the representation $\rho$, we exhibit a canonical isomorphism
between the group of $\Gal(K/\Q)$-coinvariants of $I_G X_L/I_G^2 X_L$ and the
Eisenstein ideal $\I$ modulo its square (Theorem \ref{classiso}).
Under this isomorphism,
the Frobenius element on the Galois side is identified with $U_p-1$ on the Hecke side.  This yields that the corresponding subquotient of the Galois group $Y_L$ of the maximal unramified abelian pro-$p$ extension of $L$ in which all primes above $p$ split completely is isomorphic to $\I/\I^2$ modulo the pro-$p$ subgroup generated by $U_p-1$ (Theorem \ref{Sclassiso}).
The structure of these latter subquotients can be seen to relate to the $\omega^{1-k}$-eigenspace of the quotient of the classical Iwasawa module $X_K$ by the submodule generated by an inverse limit of cup products up the cyclotomic tower.  (The precise result, Corollary \ref{specificcase}, involves a certain reciprocity map constructed out of this inverse limit of cup products in
Section \ref{cupproducts}.)
In particular, supposing $p$ does not simultaneously divide $B_k$ and $B_{p+1-k}$ for any even $k < p$, we prove that $(p,\,\cdot\,)$ is surjective if and only if $U_p-1$ generates $\I$ for each $k$ (Theorem \ref{pairgen}).  
Thus, to obtain the surjectivity in Theorem \ref{mainresult}, we verify the latter condition for all $p < 1000$.  
Note that we do not conjecture that pairing with $p$ is itself always surjective, nor that $U_p-1$ always generates $\I$.

In Section \ref{selmer}, we give a short discussion of the computation of certain Selmer groups of Hida representations taken modulo the Eisenstein ideal that is made possible by the isomorphisms of Theorems \ref{classiso} and \ref{Sclassiso}.

In forthcoming work, we will discuss an as-yet-conjectural relationship between $(\eta_i,\eta_{k-i})$ for special cyclotomic
$p$-units $\eta_i$ and odd integers $i$ (see \cite[Section 5]{mcs}), and $L_p(f,\omega^{i-1},1)$, which may be viewed as an element of $\mf{p}$, modulo $\mf{p}^2$.  Here, $f$ and $\mf{p}$ are as before (and again, for simplicity, $p^2 \nmid B_k$).  In this setting, we have $p = \eta_1$, and $L_p(f,1)$ is explicitly related to $U_p-1$ by \cite[Theorem 4.8]{kitagawa}.  Conjecture \ref{pairingconj} may be reinterpreted as the statement that the images of these $p$-adic $L$-values generate $\mf{p}/\mf{p}^2$  for each $k$.  Theorems \ref{identification} and \ref{pairgen} serve as the first pieces of theoretical evidence for this more powerful conjecture.  This conjecture suggests that the cup product pairing is only mildly degenerate, which is supported by the numerical evidence of \cite[Theorem 5.1]{mcs}.

\begin{ack}
The author thanks Sasha Goncharov, Dick Gross, Manfred Kolster, Barry Mazur, and Bill McCallum for helpful conversations and encouragement.  He also thanks Yoshitaka Hachimori and Barry Mazur for comments on and corrections to drafts of the paper.  

This material is based upon work supported by the National Science Foundation under Grant No.\ DMS-0102016 and the Natural Sciences and Engineering Research Council of Canada. This research was undertaken, in part, thanks to funding from Harvard University, the Max Planck Institute for Mathematics, and the Canada Research Chairs Program. 
\end{ack}

\section{Preliminaries} \label{prelim}

We begin by describing the ordinary Hecke algebras of Hida, which Hida defined and studied in a series of papers (e.g., \cite{hida, hida2}).  
Let $p \ge 5$ be a prime number, and let $N$ be a positive integer prime to $p$.
The ordinary Hecke algebras parameterize
ordinary modular forms of all weights and all levels $Np^r$ with $r \ge 1$. Let us very briefly recall one definition of these Hecke algebras here.

Let $M_k(\Gamma_1(Np^r))$ denote the space of modular forms over ${\bf C}$ of level $Np^r$ and weight $k \ge 2$, and let $\mf{H}_k(\Gamma_1(Np^r))$ denote the Hecke subalgebra over $\Z$ 
of $\End_{{\bf C}} M_k(\Gamma_1(Np^r))$.  We have natural maps
$$
    \mf{H}_k(\Gamma_1(Np^{r+1})) \to \mf{H}_k(\Gamma_1(Np^r)).
$$
For any commutative ring $R$ with $1$, we have a Hecke algebra
$$
    \mf{H}(N;R) = \lim_{\substack{\leftarrow \\ r \ge 1}}
    (\mf{H}_k(\Gamma_1(Np^r)) \otimes_{\Z} R)
$$
in the inverse limit, and the Hecke operators in this inverse limit are the inverse limit of the Hecke operators at finite level.
We may repeat this construction with $M_k(\Gamma_1(Np^r))$
replaced by $S_k(\Gamma_1(Np^r))$, the space of cusp forms
of weight $k$ and level $Np^r$ over ${\bf C}$.  We obtain Hecke
algebras denoted $\mf{h}_k(\Gamma_1(Np^r))$ and $\mf{h}(N;R)$.

Set
$$
  \Z_{p,N} = \lim_{\substack{\leftarrow \\ r}} \Z/p^rN\Z.
$$
If $\mc{O}$ is a $\zp$-algebra, then $\mf{H}(N;\mc{O})$ and
$\mf{h}(N;\mc{O})$ are algebras over $\mc{O}[[\Z_{p,N}^{\times}]]$
via the action of the inverse limits of diamond operators. Let
$[x] \in \mc{O}[[\Z_{p,N}^{\times}]]$ denote the group element
associated to $x \in \Z_{p,N}^{\times}$. We use $T_l$ for $l \nmid
Np$ and $U_l$ for $l \mid Np$ to denote the usual Hecke operators
for prime $l$ in these algebras.

Denote by $H^1(N;\mc{O})$ the inverse limit of the cohomology
groups $H^1(Y_1(Np^r),\mc{O})$, where $Y_1(Np^r)$ is the usual
modular curve of level $Np^r$.
The ``ordinary parts" of these cohomology groups (i.e., the
summands upon which $U_p$ acts invertibly) were studied by Ohta
in a series of works \cite{ohta-eich}-\cite{ohta-comp}, and the
rest of this section is a partial summary of his results.

\begin{remark}
  Following Ohta \cite{ohta}, we
  consider the action of $\mf{H}(N;\mc{O})$ on $H^1(N;\mc{O})$ which arises
  as the inverse limit of the actions of its Hecke operators dual under Weil
  pairings to the usual actions at each level $Np^r$.  However, we do not
  give the Hecke operators acting in this manner a separate notation
  (e.g., $T(l)^*$), since it is the only action on $H^1(N;\mc{O})$ we shall use.
  We also use the original $\Z_{p,N}^{\times}$-action of Hida 
  \cite[p.\ 552]{hida},
  which is to say that the $\Z_{p,N}^{\times}$-action
  of Ohta \cite[Section 2.1]{ohta-eich} on the Hecke algebra $\mf{H}(N;\mc{O})$ is
  a twist of ours by multiplication by $\chi^{-2}$.
  For any prime $l$ not dividing $Np$, the operator
  $T(l,l)$ of Hida 
  satisfies $[l] = l^2T(l,l)$
  in our notation.
\end{remark}

From now on, we fix a nontrivial even Dirichlet character
$\theta \colon (\Z/Np\Z)^{\times} \to \overline{\qp}^{\times}$ of
conductor $N$ or $Np$.   We let $\mc{O}$ be the ring generated over $\zp$ by the
values of $\theta$.  Let $\omega \colon \Z_{p,N}^{\times} \to
\zp^{\times}$ denote the natural Dirichlet (Teichm\"uller) character factoring
through $(\Z/p\Z)^{\times}$. Let $\chi \colon \Z_{p,N}^{\times}
\to \zp^{\times}$ denote the natural projection, and let $\kappa =
\chi\omega^{-1}$.   We will view $\kappa$, $\chi$, $\omega$, and
$\theta$ alternatively as characters $G_{\Q} \to \mc{O}^{\times}$.
For any Dirichlet character $\psi$, we let $B_{1,\psi}$ denote the associated
(first) generalized Bernoulli number.

The following conditions are employed by Ohta in his work.

\begin{hypothesis} \label{hyp} We shall make the following assumptions:
\begin{enumerate}
    \item[a.] $p \nmid \varphi(N)$, where $\varphi$ denotes the Euler-phi
    function,
    \item[b.] $\theta|_{(\Z/p\Z)^{\times}} \neq \omega|_{(\Z/p\Z)^{\times}}$,
    \item[c.] $\theta|_{(\Z/p\Z)^{\times}} \neq \omega^2|_{(\Z/p\Z)^{\times}}$,
    \item[d.] $p \mid B_{1,\theta\omega^{-1}}$,
    \item[e.] $p \nmid B_{1,\omega\theta^{-1}}$.
\end{enumerate}
\end{hypothesis}

\begin{remark}
	If $N = 1$, then Hypothesis \ref{hyp} reduces to parts $d$ and $e$.
\end{remark}

The Eisenstein ideal $\II_{\theta}(N;\mc{O})$ (resp.,
$\I_{\theta}(N;\mc{O})$) for the character $\theta$ is defined to
be the ideal of $\mf{H}(N;\mc{O})$ (resp., $\mf{h}(N;\mc{O})$)
generated by $T_l - 1 - l^{-1}[l]$  
and $[l]-\theta(l)[\kappa(l)]$
for $l \nmid
Np$, along with $U_l - 1$ for $l \mid Np$.   Note that Hypothesis \ref{hyp}a implies 
that $\mc{O}$ is unramified over $\zp$.  So, let $\mf{M}$
(resp., $\mf{m}$)
be the maximal ideal
of $\mf{H}(N;\mc{O})$ containing $\II_{\theta}(N;\mc{O})$ (resp.,
the maximal ideal of $\mf{h}(N;\mc{O})$ containing
$\I_{\theta}(N;\mc{O})$). 

We let $\mf{H} = \mf{H}(N;\mc{O})_{\mf{M}}$, the localization of
$\mf{H}(N;\mc{O})$ at $\mf{M}$, and similarly, we let $\mf{h} =
\mf{h}(N;\mc{O})_{\mf{m}}$. Let $\II$ (resp., $\I$) denote the
image of $\II_{\theta}(N;\mc{O})$ (resp., $\I_{\theta}(N;\mc{O})$)
in $\mf{H}$ (resp., $\mf{h})$.   The Hecke algebras $\mf{H}$ and $\mf{h}$ 
have $\theta$-actions of the diamond operators in $(\Z/Np\Z)^{\times}$.
We define $\Lambda_{\mf{h}}$ to be the $\zp$-subalgebra $\zp[[1+p\zp]]$ of $\mf{h}$ topologically generated by the (images of the) diamond operators in $\mf{h}$, and we may identify this with the power series ring $\zp[[T]]$ by identifying
$[1+p]-1$ and $T$.  We have an isomorphism of Hecke
$\Lambda_{\mf{h}}$-modules
$$
  \mf{h}/\I \cong \Lambda_{\mf{h}}/(g_{\theta}),
$$
where $g_{\theta}((1+p)^s-1) = L_p(\theta,1-s)$ for all $s \in
\zp$ \cite[Corollary A.2.4]{ohta-cong}.

Let $\mf{X} =H^1(N;\mc{O}) \otimes_{\mf{H}}
\mf{h}$, which is both a $G_{\bf Q}$ and an $\mf{h}$-module.  Fix a decomposition group $D_p$ and its inertia subgroup
$I_p$ at $p$ inside the absolute Galois group $G_{\Q}$ of $\Q$.
Define $\mf{X}^+$ to be the $D_p$-module $\mf{X}^{I_p}$. (Propositions 1.2.8
and 1.3.6 of \cite{ohta2} together imply that $\mf{X}^+$ is
$D_p$-stable.)  Pick a cyclic subgroup $\Delta_p$ of $I_p$ having
order $p-1$ and on which $\omega$ is injective.  We let $\mf{X}^-$ be
the $\zp$-submodule of $\mf{X}$ upon which $\Delta_p$ acts by
$\theta^{-1}\omega$.  We have a direct sum decomposition $\mf{X} = \mf{X}^+
\oplus \mf{X}^-$ of $\mf{h}$-modules (as in \cite[(5.3.5)]{ohta}).  
Ohta \cite[Theorem 3.3.2]{ohta-comp} showed that our hypotheses, particularly Hypothesis \ref{hyp}e, imply that $\mf{H}$ is Gorenstein (which also follows
from previous work of C.\ Skinner and A.\ Wiles \cite{sw}).  From this, one obtains that $\mf{X}^+$ and $\mf{X}^-$ are free $\mf{h}$-modules of rank $1$ \cite[Section 3.4]{ohta-comp}.

Let
$$
  \rho \colon G_{\Q} \to \Aut_{\mf{h}}(\mf{X})
$$
denote the Galois representation on $\mf{X}$.  The fixed field of the
kernel of $\rho$ is unramified over $\Q$ outside of the primes
dividing $Np$. We have
$$
  \det \rho(\sigma) =
  \theta(\sigma)^{-1}\chi(\sigma)[\kappa(\sigma)^{-1}] \in \mf{h}.
$$
We have four maps on $G_{\Q}$ called $a$, $b$, $c$, and $d$, with
\begin{eqnarray*}
    a(\sigma) \in \End_{\mf{h}}(\mf{X}^-), &  & b(\sigma) \in \Hom_{\mf{h}}(\mf{X}^+,\mf{X}^-),\\
    c(\sigma) \in \Hom_{\mf{h}}(\mf{X}^-,\mf{X}^+), & \mr{\ and\ } & d(\sigma) \in \End_{\mf{h}}(\mf{X}^+)
\end{eqnarray*}
for $\sigma \in G_{\Q}$, such that we may represent $\rho$ in matrix form as
$$
  \rho(\sigma) = \left(
  \begin{array}{cc}
    a(\sigma) & b(\sigma) \\
    c(\sigma) & d(\sigma)
  \end{array} \right).
$$
We remark that for $\sigma \in I_p$, we have
$$
  \rho(\sigma) = \left(
  \begin{array}{cc}
    \det \rho(\sigma) & 0 \\
    c(\sigma) & 1
  \end{array} \right),
$$
and $\rho$ is lower-triangular on an element of a
decomposition group $D_p$ containing $I_p$.

Note that we have canonical isomorphisms $\End_{\mf{h}}(\mf{X}^{\pm}) \cong
\mf{h}$ and noncanonical isomorphisms
$\Hom_{\mf{h}}(\mf{X}^{\pm},\mf{X}^{\mp}) \cong \mf{h}$ given by choices of
generators
of $\mf{X}^+$ and $\mf{X}^-$, which we fix.  With these
identifications, let $B$ and $C$ denote the
ideals of $\mf{h}$ generated by the images of $b$ and $c$,
respectively, under these isomorphisms.  
Ohta \cite[Section 3.4]{ohta-comp} used the above-mentioned freeness of $\mf{X}^{\pm}$ to show that $B = \I$ and $C = \mf{h}$ using a method developed by Harder-Pink \cite{hp} and Kurihara \cite{kurihara} for level $1$ modular forms.

Let $F$ be the fixed field of the kernel of $\theta$ on
$G_{\Q(\mu_p)}$.  Let $K$ denote the cyclotomic $\zp$-extension of
$F$.  Consider the maps
\begin{eqnarray*}
  \bar{a} \colon G_{\Q} \to (\mf{h}/\I)^{\times} & \mr{and} &
  \bar{b} \colon G_{\Q} \to B/\I B
\end{eqnarray*}
given by following $a$ and $b$ by the obvious projections.  The
fixed field of the kernel of
$$
  \bar{a} = \det \rho \pmod{\I}
$$
(see \cite[Lemma 3.3.5]{ohta2}) on $G_{\Q(\mu_p)}$ is $K$.  We have a homomorphism
$$
    \phi_B \colon G_{\Q} \to B/\I B \rtimes (\mf{h}/\I)^{\times},
$$
the semi-direct product taken with respect to the natural action,
which we can represent in matrix form as
\begin{equation*}
  \phi_B(\sigma) = \left(
    \begin{array}{cc}
      \bar{a}(\sigma) & \bar{b}(\sigma) \\
      0 & 1
    \end{array} \right)
\end{equation*}
for $\sigma \in G_{\Q}$.  Let $H$ denote the fixed field of the kernel
of $\phi_B$ on $G_F$.

Let $\Gamma = \Gal(K/F)$, let $\tgam = \Gal(K/\Q)$, let $\Lambda =
\mc{O}[[\Gamma]]$, and let $\tlam = \mc{O}[[\tgam]]$.  Let
$\tdel = \Gal(F/\Q)$, and let $\Delta =
\Gal(\Q(\mu_p)/\Q)$.  Note that $\chi$, $\kappa$, $\theta$, and
$\omega$ induce characters on $\tgam$, $\Gamma$, $\tdel$,
and $\Delta$, respectively.  We make an identification $\Lambda
\cong \mc{O}[[U]]$ for $U = \gamma-1$, where $\gamma \in \Gamma$
satisfies $\chi(\gamma) = 1+p$. 

Let $A$ be a $\zp[[\tgam]]$-module.
Given a character $\psi \colon \tdel \to \overline{\qp}^{\times}$,  
we let $A^{(\psi)}$ denote the $\psi$-eigenspace of $A$.
That is, if we let
$R_{\psi}$ denote the ring generated by the values of $\psi$ over $\zp$, then
$$
    A^{(\psi)} = A \otimes_{\zp[\tdel]} R_{\psi} 
    \cong \{ a \in A \otimes_{\zp} R_{\psi} \mid \delta a = \psi(\delta)
   a \mr{\ for\ all\ } \delta \in \tdel \},
$$
where the map $\zp[\tdel] \to R_{\psi}$ in the tensor product is the surjection
induced by $\psi$.  We remark that $A^{(\psi)}$ has the structure of a
$R_{\psi}[[\tgam]]$-module.  
In Appendix \ref{eigenspaces}, we discuss the decomposition of $\zp[\tdel]$-modules and tensor products of such modules into eigenspaces.

Let $X_K$ denote the Galois group of the maximal unramified
abelian pro-$p$ extension of $K$.  Ohta proved the following result on
the structure of its $(\omega\theta^{-1})$-eigenspace. 
Though the proof
we mention here uses the Iwasawa main conjecture, Ohta also proved it independently in order
to obtain another proof of the Iwasawa main conjecture (for $p \ge 5$).

\begin{theorem}[Ohta] \label{blemma}
    The extension $H/K$ is unramified everywhere and
    $\Gal(H/K) = X_K^{(\omega\theta^{-1})}$.
    Furthermore, $\bar{b}$ induces a canonical
    isomorphism
    $$
        X_K^{(\omega\theta^{-1})} \cong B/\I B
    $$
    of $\tlam$-modules in which the
    Galois action on $B/\I B$ is given by $\bar{a}$.
    In particular, $B/\I B$ has characteristic ideal $(g_{\theta})$
    as a $\Lambda_{\mf{h}}$-module.
\end{theorem}

\begin{proof}
    The first two sentences follow from \cite[Theorem 3.3.12]{ohta2}.
    Let $f_{\theta} \in \Lambda$ satisfy
    $$
        f_{\theta}((1+p)^s-1) = L_p(\theta,s)
    $$
    for $s \in \zp$.  The Iwasawa main conjecture, as first
    proven in \cite{mw}, states that $X_K^{(\omega\theta^{-1})}$
    has characteristic ideal $(f_{\theta}(U))$ as a $\Lambda$-module.
    Thus, to verify the final sentence, we need merely note that
    $\bar{a}(\gamma) = (1+p)[1+p]^{-1}$
    and $g_{\theta}((1+p)^s-1) = f_{\theta}((1+p)^{1-s}-1)$.
\end{proof}

\section{The Galois Representation} \label{galrep}

In this section, we study the maps
\begin{eqnarray*}
  \bar{c} \colon G_{\Q} \to C/\I C & \mr{and} &
  \bar{d} \colon G_{\Q} \to (\mf{h}/\I^2)^{\times}
\end{eqnarray*}
arising from $c$ and $d$ by projection.  The fixed field of $\bar{c}$ of $G_K$ contains a pro-$p$ abelian extension $L$ of $K$ unramified outside and totally ramified at primes above $p$, with Galois group $G$.  
The image of $d$ is contained in $1 + \I$ \cite[Lemma 3.3.5]{ohta2}.  
We show that $\bar{d}-1$ sets up a canonical isomorphism between
$(I_G X_L/I_G^2 X_L)_{\tgam}$ and $\I/\I^2$, where $X_L$ denotes the Galois group of the
maximal unramified abelian pro-$p$ extension of $L$.  If we replace $X_L$ by its
maximal quotient $Y_L$ in which all primes above $p$ split completely, then we obtain
an isomorphism with the quotient of $\I/\I^2$ by the image of the pro-$p$ group
generated by $U_p-1$.

For any algebraic extension $E$ of $\Q$, 
we let $S_E$ denote the set of primes above $p$ and any real infinite places, and we let $G_{E,S}$ denote the Galois group of the maximal unramified outside $S_E$ extension of $E$.

\begin{lemma} \label{triveigspace}
    The $\Lambda$-module $X_K^{(\theta\omega^{-1})}$ is trivial,
    and $X_K^{(\omega\theta^{-1})}$ is free of finite rank over $\mc{O}$.
\end{lemma}

\begin{proof}	Since $\theta$ is even, both
	$X_K^{(\theta\omega^{-1})}$ and $X_K^{(\omega\theta^{-1})}$ have no finite 	$\Lambda$-submodules
	(see \cite[Proposition 13.28]{washington}).  Furthermore, they both have trivial
	$\mu$-invariant by \cite{fw}, and therefore they are both free of finite
	rank over $\mc{O}$.  Hypothesis \ref{hyp}e (together with a, b, and c) 
	implies that the $p$-adic power series associated to 
	$L_p(\omega^2\theta^{-1},s)$ is a unit (see \cite[Theorem 5.11]{washington}).  	Iwasawa's Main Conjecture then forces
	$X_K^{(\theta\omega^{-1})}$ to be finite, hence trivial.
\end{proof}

Now, let us consider the homomorphism $\phi_C \colon G_{\Q} \to
C/\I C \rtimes (\mf{h}/\I)^{\times}$ that is defined by its matrix
representation
\begin{equation} \label{phiC}
  \phi_C(\sigma) = \left(
    \begin{array}{cc}
      \bar{a}(\sigma) & 0 \\
      \bar{c}(\sigma) & 1
    \end{array} \right)
\end{equation}
(or, equivalently, by the action of $G_{\Q}$ on $\mf{X}/\I \mf{X}$). Let $L$
denote the fixed field of the kernel of $\phi_C$ on $G_F$, and let
$G = \Gal(L/K)$.  We remark that $G \cong G^{(\theta\omega^{-1})}$, which
gives $G$ the structure of a $\tlam$-module.

For $E/K$, let $Z_E$ denote the Galois group of the maximal abelian pro-$p$ 
unramified outside $p$ extension of $E$. 
We shall require the following lemma.

\begin{lemma} \label{tame}
	The group $X_K^{(\omega\theta^{-1})}$ has no 
	$\mc{O}[[\Gamma']]$-submodule isomorphic to $\mc{O}(-1)$ for any
	nontrivial subgroup $\Gamma'$ of $\Gamma$.
\end{lemma}

\begin{proof}
	Let $K_n = F(\mu_{p^n})$ for $n \ge 1$, and let $A_{K_n}$ denote the
	$p$-part of the class group of $K_n$.
	We have a canonical exact sequence of $\zp[\Gal(K_n/\Q)]$-modules:
	\begin{equation} \label{h2seq}
		0 \to A_{K_n} \otimes \mu_{p^n} \to H^2(G_{K_n,S},\mu_{p^n}^{\otimes 2})
		\to \bigoplus_{v \in S_{K_n}} \mu_{p^n} \to \mu_{p^n} \to 0.
	\end{equation}
	Taking inverse limits,
	we obtain a canonical injection
	\begin{equation} \label{xkh2}
		X_K^{(\omega\theta^{-1})}(1) \hookrightarrow
		\lim_{\leftarrow} H^2(G_{K_n,S},\mu_{p^n}^{\otimes 2})^{(\omega^2
		\theta^{-1})}
	\end{equation}
	(that is an isomorphism under Hypothesis \ref{hyp}b).
	The group of $\Gal(K/K_n)$-coinvariants
	of the latter term of \eqref{xkh2} is isomorphic to the
	$(\omega^2\theta^{-1})$-eigenspace of $H^2(G_{K_n,S},\zp(2))$,
	as $G_{F,S}$ has $p$-cohomological dimension $2$.  Since 
	$\omega^2\theta^{-1}$ is even, we may replace $K_n$ with its maximal
	totally real subfield $K_n^+$, and $H^2(G_{K_n^+,S},\zp(2))$ is
	finite by a well-known result of C.\ Soul\'e \cite{soule}, as desired.
\end{proof}

Note that Lemma \ref{tame} may be stated equivalently as
saying that $f_{\theta}(U)$ is not divisible by any nontrivial divisor of 
$((U+1)^{p^n}-(1+p)^{-p^n})$ for any $n$.
We are now able to prove the following lemma on the representation $\phi_C$.

\begin{lemma} \label{clemma}
    The extension $L/K$ is unramified outside $p$ and totally
    ramified at all primes above $p$ in $K$.  The map $\bar{c}$ 
    induces an isomorphism $G \cong C/\I C$ of
    $\tlam$-modules, the action of $\tgam$ on $C/\I C$ being
    given by $\bar{a}^{-1}$.  We also have $C/\I C \cong
    \Lambda_{\mf{h}}/(g_{\theta}(T))$ as
    $\Lambda_{\mf{h}}$-modules.
\end{lemma}

\begin{proof}
    The statement that $G \cong C/\I C$ follows
    as in \cite[p.\ 298]{ohta}, the statement on the action
    is obvious, and the final statement is also obvious as $C = \mf{h}$.
    Note that we therefore have an isomorphism
    $G \cong \Lambda/(f^{\circ}_{\theta}(U))$ of $\Lambda$-modules,
    where $f_{\theta}^{\circ}(U) = f_{\theta}((U+1)^{-1}-1)$. 

    As for the first statement, we note that if $v$ is a prime of $K$ lying 
    above $N$, the inertia group $I_v$ of $v$ in $G$ cannot
    be finite and nontrivial, as the $\Lambda$-submodule generated by $I_v$
    would then be finite (as the set of primes of $K$ above $N$ is finite) 
    and therefore trivial.  On the other hand, any ramification at $v$ is tame, 
    so if $I_v$ is infinite, it must be
    isomorphic to $\zp(1)$ as a module over $\zp[[\Gamma']]$ for some
    subgroup $\Gamma' \le \Gamma$ of finite index.  
    This would imply that
     $f_{\theta}^{\circ}(U)$ shares a nontrivial common divisor with
    $(U+1)^{p^n}-(1+p)^{p^n}$ for some $n$, which is not the case by Lemma \ref{tame}.  
    Hence, we see that $L/K$ is unramified outside $p$.
    That $L/K$ is totally ramified at primes above $p$ now follows from
    Lemma \ref{triveigspace} and the fact that $G$ has a 
    $(\theta\omega^{-1})$-action of $\tdel$.
\end{proof}

For any algebraic extension $E$ of $K$, let $X_E$ denote the
Galois group of the maximal unramified abelian pro-$p$ extension
of $E$, and let $Y_E$ denote the maximal quotient of $X_E$ in
which all primes above $p$ split completely. Let $I_G$ denote the
augmentation ideal in $\zp[[G]]$. We have a canonical isomorphism
$G \cong I_G/I_G^2$ of $\zp[[\tgam]]$-modules by $\sigma \mapsto \sigma-1$
(and hence of $\tlam$-modules, as $G$ has the structure of an $\mc{O}$-module).

From now on, we assume the following strengthening of Hypothesis \ref{hyp}b. 

\begin{hypothesis} \label{square}
	We have that
	$\theta^2|_{(\Z/p\Z)^{\times}} \neq \omega^2|_{(\Z/p\Z)^{\times}}$.
\end{hypothesis}

\begin{remark}
	In fact, it would suffice in what follows
	to assume that $\theta^2$ and $\omega^2$ do not agree on $D_p$.
	If $N = 1$, then Hypothesis \ref{square} holds automatically.
\end{remark}  

\begin{lemma} \label{samemodules}
    The canonical maps in the following commutative diagram are all $\tlam$-isomorphisms:
    $$ \SelectTips{eu}{}
        \xymatrix{
        (X_L/I_G X_L)^{(\omega\theta^{-1})} \ar[r]^{\qquad \sim} \ar[d]^{\wr}
        &  X_K^{(\omega\theta^{-1})} \ar[d]^{\wr} \\
        (Y_L/I_G Y_L)^{(\omega\theta^{-1})} \ar[r]^{\qquad \sim} &
         Y_K^{(\omega\theta^{-1})}.}
    $$
\end{lemma}

\begin{proof}
    	For $v \in S_K$, let $G_v$ denote the decomposition group at $v$ in $G$,
	let $D'_v$ denote the decomposition group at $v$ in $X_L/I_G X_L$, and let $D_v$ 	
	denote the decomposition group at $v$ in $X_K$.
    	Since, by Lemma \ref{clemma}, we know that $L/K$ is unramified outside $p$
    	and has trivial maximal unramified subextension, we have a commutative diagram
    	with exact rows and columns: 
   	 $$
        	\SelectTips{eu}{} \xymatrix@R=15pt@C=15pt{
        	\bigoplus_{v \in S_K} D'_v \ar[r] \ar[d] & \bigoplus_{v \in S_K} D_v \ar[d] & \\
        	X_L/I_G X_L \ar[r] \ar[d] &
        	X_K \ar[r] \ar[d] & 0 \\
        	Y_L/I_G Y_L \ar[r] \ar[d] & Y_K \ar[r] \ar[d] & 0 \\
       	0 & 0 & }
    	$$

    Note that $\Delta$ does not permute the primes in $S_K$.  
    Furthermore, each $D_v$ has
    a trivial action of $\Delta$, as it is the Galois group of an unramified local extension
    and the restriction of $\omega\theta^{-1}$ to $\Delta$ is not the identity 
    by Hypothesis \ref{hyp}b.  Therefore, the direct sum of the $D_v$ has trivial 
    $(\omega\theta^{-1})$-eigenspace.  The same holds for the direct sum of the
    $D'_v$.  
    
	It now suffices to show that the $(\omega\theta^{-1})$-eigenspace of the
	kernel $\mc{K}$ of $Y_L/I_G Y_L \to Y_K$ is trivial.  For this, consider the following
	commutative diagram with exact columns and exact rows aside from the first, which
	is a complex:
	$$
		\SelectTips{eu}{} \xymatrix@R=15pt@C=15pt{
		& & 0 \ar[d] & \\
		\bigoplus_{v \in S_K} H_2(G_v,\zp) \ar[d] \ar[r]^{\ \ \quad j} & 
		H_2(G,\zp) \ar[r] \ar[d] & \mc{K} \ar[d] 	&\\
		 \bigl( \bigoplus_{w \in S_L} H_1(G_{L_w},\zp) \bigr)_G 
		\ar[r] \ar[d] & (Z_L)_G \ar[r] \ar[d] & (Y_L)_G \ar[r] \ar[d] & 0\\
		\bigoplus_{v \in S_K} H_1(G_{K_v},\zp) \ar[r] \ar[d] &
		Z_K \ar[r] \ar[d] & Y_K \ar[r] \ar[d] & 0 \\
		\bigoplus_{v \in S_K} G_v \ar[r] \ar[d] & G \ar[r] \ar[d] & 0\\ 
		0 & 0 &}
	$$
	(Note that $H_2(G,\zp)$ 
	is isomorphic to the wedge product over $\zp$ of two copies of $G$.)  Since
	$G_v = G$ for each $v \in S_K$, we have that the corestriction map $j$ is
	surjective.	By a diagram chase, it remains only to show that the 
	$(\omega\theta^{-1})$-eigenspace of $\bigoplus_{v \in S_K} G_v$ is trivial.  
	For this,
	we remark as before that its summands are not permuted by $\Delta$.  Since 
	$G \cong G^{(\theta\omega^{-1})}$, Hypothesis \ref{square} implies that the 	$\zp[\Delta]$-eigenspace of $G_v$ determined by the restriction of 	$\omega\theta^{-1}$ is trivial.  
\end{proof}

We now study the homomorphism on $G_{\Q}$ to formal matrices defined as
$$
    \phi_D(\sigma) = \left(
    \begin{array}{cc}
      \bar{a}(\sigma) & \bar{b}(\sigma) \\
      \bar{c}(\sigma) & \bar{d}(\sigma)
    \end{array} \right)
$$
for $\sigma \in G_{\Q}$.  That is, consider the set
$$
	\mc{N} = \left\{ \left( \begin{array}{cc}
	\alpha & \beta \\ 
	\gamma & \delta
	\end{array} \right) \mid
	\alpha \in (\mf{h}/\I)^{\times},\ \beta \in \I/\I^2,\ 
	\gamma \in \mf{h}/\I, \text{ and } 
	\delta \in (\mf{h}/\I^2)^{\times} \right\},
$$ 
which we observe forms a group under the usual multiplication of matrices.  More specifically, for $\delta \in (\mf{h}/\I^2)^{\times}$, we let $\bar{\delta}$ denote its image in $(\mf{h}/\I)^{\times}$, and we view $\mf{h}/\I$ and $\I/\I^2$ as $\mf{h}/\I$-modules under multiplication.  Our multiplication is:
$$
	\left( \begin{array}{cc}
	\alpha & \beta \\ 
	\gamma & \delta
	\end{array} \right) \cdot
	\left( \begin{array}{cc}
	\alpha' & \beta' \\ 
	\gamma' & \delta'
	\end{array} \right) = 
	\left( \begin{array}{cc}
	\alpha\alpha' & \alpha\beta'+\overline{\delta'}\beta \\ 
	\alpha'\gamma + \bar{\delta}\gamma' & \gamma\beta'+\delta\delta'
	\end{array} \right).
$$
Then $\phi_D$ is the homomorphism given by following $\rho$ with the projection from the image of $\rho$ in $\mr{GL}_2(\mf{h})$ to $\mc{N}$, noting again
here that $B = \I$.

Let $M$ denote the fixed field of the
kernel of $\phi_D$ on $G_F$.  Although we do not need it, we note
that we may use $\phi_D$ to form a new homomorphism $\psi_D$,
represented in matrix form as
$$
  \psi_D(\sigma) = \left(
  \begin{array}{ccc}
    1 & \bar{c}(\sigma) & \bar{d}(\sigma)-1 \\
    0 & \bar{a}(\sigma) & \bar{b}(\sigma) \\
    0 & 0 & 1
  \end{array} \right)
$$
for $\sigma \in G_{\Q}$, which perhaps gives one a better feeling
for the structure of the Galois group $\Gal(M/\Q)$.  In particular, we may
view $M$ as a Heisenberg extension of $K$.  We have a field
diagram:
$$
	\SelectTips{eu}{}
    
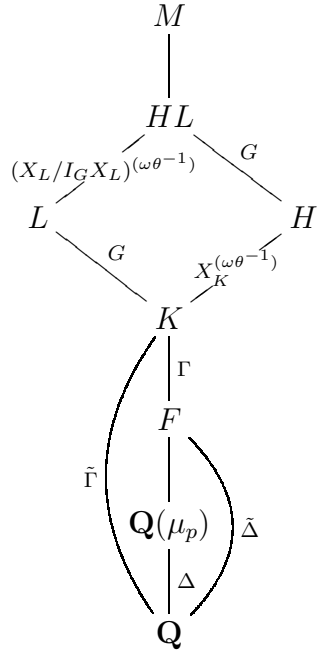
\xymatrix{
    & M \ar@{-}[d] \\
    & HL \ar@{-}[ld]|{(X_L/I_G X_L)^{(\omega\theta^{-1})}} \ar@{-}[rd]^G &&
    \\ 
    L \ar@{-}[dr]^{G} && H \ar@{-}[ld]|-{X_K^{(\omega\theta^{-1})}} \\ 
    &K 
    \ar@{-}[d]^{\Gamma} \ar@/_2pc/@{-}[ddd]_{\tilde{\Gamma}} \ar@{-}[d]
    & \\
    &F \ar@/^2pc/@{-}[dd]^{\tdel} \ar@{-}[d] && \\
    &\Q(\mu_p) \ar@{-}[d]^{\Delta}\\
    &\Q
    }
$$
We remark that $G \cong C/\I C$ has an $\bar{a}^{-1}$-action of $\tgam$, hence a $(\theta\omega^{-1})$-action of $\tdel$, as opposed
to the $\bar{a}$-action of $\tgam$ (resp., $(\omega\theta^{-1})$-action of
$\tdel$) on $\Gal(H/K) \cong B/\I B$.
Since $M/H$ and $M/L$ are abelian, $\Gal(M/HL)$ too has an action of $\tgam$, but this action is seen to be trivial from the matrix representation.

\begin{theorem} \label{classiso}
  The extension $M/HL$ is unramified everywhere, and
  $$
    \Gal(M/HL) \cong (I_G X_L/I_G^2 X_L)_{\tgam}.
  $$
  The map $\bar{d} - 1$ induces a canonical isomorphism
  of $\zp$-modules
  $$
    (I_G X_L/I_G^2 X_L)_{\tgam} \xrightarrow{\sim} \I/\I^2.
  $$
\end{theorem}

\begin{proof}
  Note the following equality in $\I/\I^2$:
  $$
    \bar{d}(\sigma\tau)-1 =
    \bar{c}(\sigma)\bar{b}(\tau)+\bar{d}(\sigma)\bar{d}(\tau)-1
    = \bar{c}(\sigma)\bar{b}(\tau)+(\bar{d}(\sigma)-1)+(\bar{d}(\tau)-1)
  $$
  for $\sigma, \tau \in G_{\Q}$.
  One easily computes that, for $\sigma \in G_H$ and $\tau \in
  G_L$, we have
  \begin{equation*} \label{bcd}
    \bar{d}([\sigma,\tau]) - 1 = \bar{c}(\sigma)\bar{b}(\tau).
  \end{equation*}
  Since $B = \I$, $C = \mf{h}$, and $\bar{b}$ and $\bar{c}$ are
  surjective on $G_L$ and $G_H$, respectively, the group
  $P = \Gal(M/HL)$ is seen to be equal to the commutator
  subgroup of $\Gal(M/K)$, and $\bar{d}-1$ induces
  a canonical isomorphism
  \begin{equation} \label{eisiso}
    P \xrightarrow{\sim} \I/\I^2.
  \end{equation}

  Now, note that $M$ is an abelian pro-$p$ extension of
  $HL$ which is unramified outside the primes dividing $N$.  
  Furthermore, the completion of $HL$ at a prime $w$ above $N$ is
  the unique unramified $\zp$-extension of the completion of $F$
  at this prime.  If the inertia group
  $P_w$ of $w$ in $P$ is infinite, then, as the ramification at
  $w$ is tame, it is isomorphic to $\zp(1)$ as a module over
  $\zp[[\Gamma']]$ for some subgroup $\Gamma' \le \Gamma$ of finite index.
  However, this is impossible as $P$ has a trivial
  $\tgam$-action.  Hence, $P_w$ must be finite.  
  However, since $P$
  is isomorphic to $X_K^{(\omega\theta^{-1})}$ as a pro-$p$ group
  by Theorem \ref{blemma} and
  \eqref{eisiso}, the second statement of Lemma \ref{triveigspace} forces $P_w$ to be   	zero.  Hence, $M/HL$ is everywhere unramified.

  We have now exhibited a surjection
  \begin{equation} \label{surjection}
    (I_G X_L/I_G^2 X_L)_{\tgam} \twoheadrightarrow P.
  \end{equation}
  Combining \eq{eisiso} and \eq{surjection} and noting Lemma \ref{triveig}, we have a 	natural
  commutative diagram of $\zp[[\tgam]]$-modules:
  \begin{equation} \label{diagram}
    \SelectTips{eu}{}
    \xymatrix{
    I_G/I_G^2 \otimes_{\mc{O}} (X_L/I_G X_L)^{(\omega\theta^{-1})} \ar[d]^{\wr} \ar@{>>}[r] &
    (I_G X_L/I_G^2 X_L)_{\tgam} \ar@{>>}[d] \\
    C/\I C \otimes_{\mc{O}} B/\I B  \ar@{>>}[r] & \I/\I^2.
    }
  \end{equation}

  In fact, the lower horizontal map in \eq{diagram} factors through
  the tensor product over $\mf{h}$.
  Of course, since $C = \mf{h}$ and $B = \I$, we have that the
  map
  \begin{equation} \label{tensorisom}
    C/\I C \otimes_{\mf{h}} B/\I B \to \I/\I^2
  \end{equation}
  is an isomorphism.  As the $G_{\Q}$-action on $B/\I B$ is given by $\bar{a}$ 
  and on $C/\I C$ by
  $\bar{a}^{-1}$, it is easy to see that
  $$
    C/\I C \otimes_{\mf{h}} B/\I B \cong (C/\I C \otimes_{\mc{O}}
    B/\I B)_{\tgam}.
  $$
  Thus, the right vertical arrow in the diagram \eq{diagram} is
  also an isomorphism.
\end{proof}

\begin{theorem} \label{Sclassiso}
  Let $\mc{V}$ denote the $\zp$-span of $U_p-1$ in $\I$.
  Then the map $\bar{d}-1$ induces an isomorphism
  $$
    (I_G Y_L/I_G^2 Y_L)_{\tgam} \to \I/(\mc{V}+\I^2)
  $$
  of $\zp$-modules.
\end{theorem}

\begin{proof}
  We know from Theorem \ref{classiso} that $X_L/I_G^2 X_L$ has a
  $\zp[[\Gal(L/\Q)]]$-quotient $Z = \Gal(M/L)$ and, therefore, also a quotient
  $Z'$, which fit into a commutative diagram with exact rows
	\begin{equation} \label{1stdiagram} \SelectTips{eu}{}
    \xymatrix{
    0 \ar[r] & (I_G X_L/I_G^2 X_L)_{\tgam} \ar[r] \ar@{->>}[d] &
    Z \ar[r] \ar@{->>}[d] & (X_L/I_G X_L)^{(\omega\theta^{-1})}
    \ar[r] \ar[d]^{\wr} & 0 \\
    0 \ar[r] & (I_G Y_L/I_G^2 Y_L)_{\tgam} \ar[r] &
    Z' \ar[r] & (Y_L/I_G Y_L)^{(\omega\theta^{-1})} \ar[r] & 0.
    }
  \end{equation}
  Let $\mf{p}$ be a prime above $p$ in $K$.  We know that $L/K$ is totally
  ramified at $\mf{p}$ by Lemma \ref{clemma}, so
  there is a unique prime $\mf{P}$ above $\mf{p}$ in $L$.
  Since $X_L/I_G^2 X_L$ is abelian, there is a unique decomposition
  group $D_{\mf{p}}$ at $\mf{P}$ in this Galois group.
  By Lemma \ref{samemodules}, the image $\bar{D}_{\mf{p}}$
  of $D_{\mf{p}}$ in $Z$ is
  contained in $(I_G X_L/I_G^2 X_L)_{\tgam}$.  Since $\tgam$ acts
  trivially on the latter module, we have that $\bar{D} = \bar{D}_{\mf{p}}$ is
  the unique decomposition group above $p$ in $Z$.

  We claim that $\bar{D}$ is equal to the kernel of the surjection
  $$
    (I_G X_L/I_G^2 X_L)_{\tgam} \twoheadrightarrow (I_G Y_L/I_G^2 Y_L)_{\tgam}.
  $$
  To see this, note that the kernel of the map $X_L/I_G^2 X_L \to
  Y_L/I_G^2 Y_L$ is the product of the $D_{\mf{p}}$, and hence the
  kernel of $I_G X_L/I_G^2 X_L \to I_G Y_L/I_G^2 Y_L$ is the
  intersection of this product with $I_G X_L/I_G^2 X_L$.  Now we know that the
  projection of this product to $Z$ is simply $\bar{D}$, and this is
  the entire kernel of the map $Z \to Z'$ by \eqref{1stdiagram} and the Snake Lemma applied to the
  diagram
  $$ \SelectTips{eu}{}
    \xymatrix{
    0 \ar[r] & I_G(I_G+I_{\Gamma}) X_L/I_G^2X_L
     \ar[r] \ar@{->>}[d] &
    I_G X_L/I_G^2 X_L \ar[r] \ar@{->>}[d] & (I_G X_L/I_G^2 X_L)_{\tgam}
    \ar[r] \ar@{->>}[d] & 0 \\
    0 \ar[r] & I_G(I_G+I_{\Gamma}) Y_L/I_G^2Y_L \ar[r] &
    I_G Y_L/I_G^2 Y_L \ar[r] &  (I_G Y_L/I_G^2 Y_L)_{\tgam} \ar[r] & 0.
    }
  $$
  Hence,
  $$
    (I_G Y_L/I_G^2 Y_L)_{\tgam} \cong (I_G X_L/I_G^2
    X_L)_{\tgam}/\bar{D}.
  $$

  On the other hand, the map $d$ is a homomorphism on our fixed decomposition
  group $D_p$ in $G_{\Q}$ which
  has image $U_p$ on the inverse of a Frobenius at $p$ 
  \cite[Theorem 3.4.2]{ohta-cong} (using Hypothesis \ref{hyp}c).
  Since $[F:\Q]$ is prime to $p$ and $L/F$ is totally ramified at
  primes above $p$, a prime-to-$p$ power of this Frobenius is
  contained in $D_p \cap G_L$.  Also, we have seen that the image
  of $D_p \cap G_L$ in $Z$ is simply $\bar{D}$.
  By Theorem \ref{classiso}, we obtain immediately that
  $$
    (I_G X_L/I_G^2 X_L)_{\tgam}/\bar{D} \cong \I/(\mc{V}+\I^2),
  $$
  finishing the proof.
\end{proof}

The following example illustrates what occurs in a case in which the class
group of $F$ has $p$-rank $1$.

\begin{example}
	Consider the case of the smallest irregular prime, $p = 37$, which divides
	the Bernoulli number $B_{32}$.   Let $N = 1$ and $\theta = \omega^{32}$.  
	Here, $X_K$ is its own $\omega^{-31}$-eigenspace, and
	$$
		\Gal(H/K) \cong X_K \cong Y_K \cong \Z_{37}(1-k)
	$$ 
	for some $k \equiv 32 \bmod 36$ (inside $\Z_{37} \times \Z/36\Z$).  
	We have $G = \Gal(L/K) \cong \Z_{37}(k-1)$.  Both $H$ and $L$ are 
	generated by cyclotomic
	$37$-units.  For instance, the degree $37$ subextension of $L/K$ is generated
	by the $37$th root of
	$$
		\eta = \prod_{i=1}^{36} (1-\zeta_{37}^i)^{i^{30}},
	$$
	where $\zeta_{37}$ is a primitive $37$th root of $1$.
	
	Theorem \ref{classiso} implies
	$$
		G \otimes_{\Z_{37}} X_K \cong I_G/I_G^2 \otimes_{\Z_{37}} X_L/I_GX_L
		\cong I_G X_L/I_G^2 X_L \cong \Z_{37}
	$$
	as $\tlam$-modules, and the Galois group $I_GX_L/I_G^2 X_L$
	is the commutator
	subgroup of the Galois group of the extension it defines over $K$.
	It turns out that $U_{37}-1$ generates $\I$, so $I_G X_L/I_G^2 X_L$
	is generated by any Frobenius at $37$.  In other words, 
	$I_G Y_L/I_G^2 Y_L = 0$ and $M = HL$.  As we shall see in Theorem 
	\ref{pairgen}, 
	this also forces the nontriviality of 
	a cup product pairing at the level of $F$ on the $37$-units $37$ and
	$\eta$.
\end{example}

\section{Cup products} \label{cupproducts}

In this section, we shall work in a somewhat more abstract setting and construct
what we refer to as the $S$-reciprocity map for an abelian $p$-power Kummer extension of number fields, which specializes to local reciprocity maps at primes in a fixed set $S$ containing the primes above $p$ and those which ramify.  The $S$-reciprocity map is constructed out of cup products in Galois cohomology with ramification restricted to $S$ and coefficients in $p$-power roots of unity.  It allows us to give a description of the second graded quotient in the augmentation filtration of the Galois group of the maximal unramified abelian pro-$p$ completely split above $p$ extension of an $S$-ramified Kummer extension of a cyclotomic $\zp$-extension of a number field.

In this section, let $p$ be any prime number, and, for the time being, let $F$ be any number field containing the $p^m$th roots of unity for some $m \ge 1$.  Let $S$ be any finite set of primes of $F$ containing those above $p$
and any real places.  
For any finite extension $E$ of $F$, let $S_E$ denote the set of primes of $E$ 
consisting of the primes above those in $S$, 
let $A_{E,S}$ denote the $p$-part of the $S_{E}$-class group
of $E$ (and similarly for any extension of $F$). 
We let $G_{E,S}$ denote the Galois group of the maximal extension of
$E$ unramified outside $S_{E}$.  Write $\mc{O}_E$ and $\mc{O}_{E,S}$ for the ring of integers and $S_E$-integers of $E$, respectively.
Write $\mr{Br}_S(E)$ for the $S_{E}$-part of the Brauer
group of $E$.
In most of the rest of the notation, the set $S$ will be understood.
We write $\mc{E}_E$ for the group of $S_E$-units $\mc{O}_{E,S}^{\times}$ in $E$.

We recall the following exact sequences for $G_{F,S}$-cohomology groups:
\begin{equation} \label{h1seq}
	0 \to \mc{E}_F/\mc{E}_F^{p^m} \to H^1(G_{F,S},\mu_{p^m})
	\to A_{F,S}[p^m] \to 0
\end{equation}
and
\begin{equation} \label{brauer}
    0 \to A_{F,S}/p^m \to H^2(G_{F,S},\mu_{p^m}) \to
    \mr{Br}_S(F)[p^m] \to 0.
\end{equation}
By Kummer theory and \eqref{h1seq}, the cohomology group $H^1(G_{F,S},\mu_{p^m})$ is isomorphic to the subgroup $\mc{B}_{m,F}$ of $F^{\times}/F^{\times p^m}$ consisting of reductions of elements $a \in F^{\times}$
with fractional ideal $a\mc{O}_{F,S}$ a $p^m$th power.
We will consider the pairing
$$
    \langle\, \cdot \, ,\, \cdot \, \rangle_{m,F} \colon \mc{B}_{m,F} \times
    \mc{B}_{m,F}
    \to H^2(G_{F,S},\mu_{p^m}^{\otimes 2})
$$
obtained by cup product on $H^1(G_{F,S},\mu_{p^m})$.  

Let $\Upsilon$ be any subgroup of $\mc{B}_{m,F}$
that is free over $\Z/p^m\Z$.
Let $E$ denote the field defined over $F$ by $p^m$th roots of lifts of 
the elements of $\Upsilon$, and let $Q = \Gal(E/F)$.  
Kummer theory provides a
canonical isomorphism
\begin{eqnarray} \label{Kummer}
    Q \xrightarrow{\sim} \Hom(\Upsilon,\mu_{p^m}),&&
    \sigma \mapsto \bigl(a \mapsto \frac{\sigma(\alpha)}{\alpha}\bigr),
\end{eqnarray}
where for $a \in \Upsilon$, we use $\alpha$ to denote a $p^m$th root of
a lift of $a$.
Let $I_{Q}$ denote the augmentation ideal of $\zp[Q]$.  
We have that $Q \cong I_{Q}/I_{Q}^2$ in the usual manner, 
$\sigma \mapsto \sigma-1$.
Tensoring \eq{Kummer} with $H^2(G_{F,S},\mu_{p^m})$ 
and applying the inverse of the resulting map to the 
homomorphism
$$
    a \mapsto \langle a,b \rangle_{m,F}
$$
for $a \in \Upsilon$ and a fixed 
$b \in \mc{E}_F/\mc{E}_F^{p^m}$ (for instance), 
we obtain an element of 
$$
	I_{Q}/I_{Q}^2 \otimes H^2(G_{F,S},\mu_{p^m}).
$$
By varying $b$ and lifting to $\mc{E}_F$,
we have therefore defined a canonical homomorphism
$$
  \Psi_{m,E/F} \colon \mc{E}_F
  \to I_{Q}/I_{Q}^2 \otimes H^2(G_{F,S},\mu_{p^m}),
$$
which one might refer to as the $m$th {\em $S$-reciprocity map} for the extension $E/F$, 
as local reciprocity maps can be constructed in an analogous fashion from norm 
residue symbols.  (Note that we need not pass through the subgroup $\Upsilon$ or assume that $F$ contains $\mu_{p^m}$ in our construction, but we do so for comparison
with the cup product pairing.)

Let $\mc{U}_{E/F} = \mc{E}_F \cap N_{E/F} E^{\times}$. 
We have a $p$-group
$$
    \mc{P}_{m,E/F} = \Psi_{m,E/F}(\mc{U}_{E/F}),
$$
and we remark that
$$
	\mc{P}_{m,E/F} \subseteq I_Q/I_Q^2 \otimes \mc{A}_{m,E/F},
$$
where $\mc{A}_{m,E/F}$ is defined to be the image of the norm map
$N_{E/F} \colon A_{E,S} \to A_{F,S}/p^m$.  
(We can take $\mc{A}_{m,E/F}$ here
instead of $H^2(G_{F,S},\mu_{p^m})$, or even $A_{F,S}/p^m$, 
by the fact that global norms are local
norms everywhere and \cite[Theorem 2.4]{mcs}.
Note that $(A_{F,S}/p^m)/\mc{A}_{m,E/F}$ is canonically isomorphic to a quotient of
$Q$.)
The following is immediate from the definition of $\Psi_{m,E/F}$.

\begin{lemma} \label{eltdef}
	Let $b \in \mc{U}_{E/F}$.
  For $N \ge 0$ and for $\sigma_i \in Q$ and
  $\mf{a}_i \in \mc{A}_{m,E/F}$ with $1 \le i \le N$, we have
  $$
    \Psi_{m,E/F}(b) = \sum_{i=1}^N (\sigma_i-1) \otimes
    \mf{a}_i
  $$
  if and only if
  $$
    \langle a,b \rangle_{m,F} =
    \sum_{i=1}^N \mf{a}_i \otimes g_{\sigma_i}(a)
  $$
  for all $a \in \Upsilon$, where, for $\sigma \in Q$, we have used
  $g_{\sigma}$ to denote the image of $\sigma$ under the map
  in \eqref{Kummer}.
\end{lemma}

Let $\mc{D}_{m,E/F}$ be the kernel of $A_{E,S} \to A_{F,S}/p^m$.   We define
$$
	\mc{Q}_{m,E/F} = I_{Q} A_{E,S}/I_Q \mc{D}_{m,E/F},
$$
which is a quotient of $I_Q A_{E,S}/(p^m+I_Q)I_Q A_{E,S}$.

\begin{theorem} \label{finclassgpiso}
	There is a canonical isomorphism
  $$
    (I_{Q}/I_{Q}^2 \otimes \mc{A}_{m,E/F})/\mc{P}_{m,E/F} \xrightarrow{\sim} 
    \mc{Q}_{m,E/F}.
  $$
\end{theorem}

\begin{proof}
  We have a natural map
  $$
  	f \colon I_{Q}/I_{Q}^2 \otimes A_{E,S}/I_{Q}A_{E,S}
  	\to I_{Q} A_{E,S}/(p^m+I_{Q})I_{Q} A_{E,S}
  $$
  with 
  $$
    f((\sigma - 1) \otimes \bar{\mf{A}}) =
    (\sigma - 1)\mf{A}  \pmod{(p^m + I_{Q})I_{Q} A_{E,S}}.
  $$
  for $\sigma \in Q$ and $\mf{A} \in A_{E,S}$ with image $\bar{\mf{A}}$
  in $A_{E,S}/I_{Q} A_{E,S}$.
  This descends to a homomorphism
  $$
  	h \colon I_{Q}/I_{Q}^2 \otimes \mc{A}_{m,E/F}
  	\to \mc{Q}_{m,E/F}.
  $$
  We must show that $h$ has kernel $\mc{P}_{m,E/F}$.  
	
  Let $J_{E,S}$ denote the $S_{E}$-ideal group of $E$.
  Let $b \in \mc{U}_{E/F}$,
  and write $b = N_{E/F} y$ for some $y \in E^{\times}$.
  The ideal of $\mc{O}_{E,S}$ generated by $y$ may be written in the form
  $$
    y\mc{O}_{E,S} = \prod_{i=1}^t
    \mf{A}_i^{1-\sigma_i} 
  $$
  with $t \ge 0$, $\mf{A}_i \in J_{E,S}$, and
  $\sigma_i \in Q$.  For $a \in \Upsilon$, let
  $E_a = F(\alpha)$, where $\alpha^{p^n} = a$.
  The norm of $y$ to $E_a$ generates
  $$
    \prod_{i=1}^t (N_{E/E_a}\mf{A}_i)^{1-\sigma_i|_{E_a}}.
  $$
  Then \cite[Theorem 2.4]{mcs} implies that
  \begin{equation} \label{repeatedformula}
    \langle a,b \rangle_{m,F} =
    \sum_{i=1}^t N_{E/F}\mf{A}_i \otimes
    g_{\sigma_i}(a)
  \end{equation}
  (see also \cite[Theorem 4.3]{me-massey}).  By
  Lemma \ref{eltdef}, we conclude that $h(\Psi_{m,E/F}(b))$ is equal to the
  image of
  $$
    \sum_{i=1}^N (\sigma_i-1)\mf{A}_i
  $$
  in $\mc{Q}_{m,E/F}$,
  which is $0$ by the principality of $y$.  Thus, $\mc{P}_{m,E/F}$ is
  contained in the kernel of $h$.

  On the other hand, if
  $$
    \sum_{i=1}^N (\sigma_i-1)\mf{A}_i  $$
  has trivial image in $\mc{Q}_{m,E/F}$,
  then there exist $\mf{B}_j \in J_{E,S}$ and $\tau_j \in Q$ for some $1 \le j \le M$ and $M \ge 0$
  with $N_{E/F} \mf{B}_j$ having trivial class in $A_{F,S}/p^m$
  and such that
  $$
    \mf{C} = \prod_{i=1}^N \mf{A}_i^{1-\sigma_i}\prod_{j=1}^M
    \mf{B}_j^{1-\tau_j}
  $$
  is principal.
  Let $y \in E^{\times}$ be a generator of $\mf{C}$, and set $b = N_{E/F} y \in
  \mc{E}_F$.
  As before, \eq{repeatedformula} holds for any $a \in \Upsilon$.  Again by Lemma
  \ref{eltdef}, we
  conclude that
  $$
	\Psi_{m,E/F}(b) = \sum_{i=1}^N (\sigma_i-1) \otimes N_{E/F} \mf{A}_i.
  $$
  Hence, the kernel of $h$ equals $\mc{P}_{m,E/F}$.
\end{proof}

We pass to the infinite level.  
Now let $F$ be any number field and $S = S_F$ a set of primes containing those above
$p$ and any real places.  Let $K_n = F(\mu_{p^n})$ for all $n \ge 1$, and
let $K = F(\mu_{p^{\infty}})$.  
Let $\mc{B}_{m,K}$ denote the direct limit of the $\mc{B}_{m,K_n}$.
The pairings $\langle \, \cdot \, , \,\cdot \, \rangle_{m,K_n}$ for fixed $m$ and increasing $n$ induce, via the usual compatibility of cup products under corestriction, pairings in the limit
$$
	\mc{B}_{m,K} \times
	\lim_{\substack{\leftarrow \\ n}} \mc{E}_{K_n}/\mc{E}_{K_n}^{p^m}
	\to \lim_{\substack{\leftarrow\\ n}}\, H^2(G_{K_n,S},\mu_{p^m}^{\otimes 2}),
$$
with the inverse limits being taken with respect to norm and corestriction maps,
respectively. 
Let 
\begin{equation} \label{HK}
	\mc{H}_K = \lim_{\leftarrow} H^2(G_{K_n,S},\zp(1)),
\end{equation}
with the inverse limit taken with respect to corestriction maps.
Taking the inverse limit over $m$, we obtain a pairing
$$
	\langle \,\cdot\, , \,\cdot\, \rangle_K \colon \mc{B}_K \times \mc{U}_K \to 
	\mc{H}_K(1),
$$
where $\mc{B}_K$ denotes the inverse limit of the $\mc{B}_{m,K}$ (which contains
the $p$-completion of the $p$-units in $K$) and
$$
  \mc{U}_K = \lim_{\leftarrow}\, (K_n^{\times} \otimes \zp)
       \cong \lim_{\leftarrow}\, (\mc{E}_{K_n} \otimes \zp)
	\cong \lim_{\substack{\leftarrow\\m, n}} \mc{E}_{K_n}/\mc{E}_{K_n}^{p^m},
$$
the inverse limits over $n$ being taken with respect to norm maps.
Note that $\mc{U}_K$ is simply the group of universal norms of $p$-units, independent of $S$, as all primes not over $p$ in $K$ have infinite decomposition groups but trivial inertia groups over $F$.

Let $L$ be the Kummer extension of $K$ defined by all $p$-power roots of 
a pro-$p$ subgroup $\Upsilon$ of $\mc{B}_K$, which is defined even though
$\mc{B}_K$ is only contained in the $p$-completion of $K^{\times}$.
We remark that $G = \Gal(L/K)$ is necessarily 
$\zp$-torsion free.  
Since
$$
	\Hom(\Upsilon,\mc{H}_K(1)) \cong G \otimes_{\zp} \mc{H}_K,
$$
we have a natural map
\begin{equation} \label{PsiLK}
	\Psi_{L/K} \colon \mc{U}_K \to I_G/I_G^2 \otimes_{\zp} \mc{H}_K
\end{equation}
induced by
$$
	b \mapsto (a \mapsto \langle a,b \rangle_K)
$$
for $b \in \mc{U}_K$ and $a \in \Upsilon$.  We call $\Psi_{L/K}$ the $S$-reciprocity
map for $L/K$, and we remark that it satisfies the obvious analogue of Lemma 
\ref{eltdef}.
For each $n \ge 1$, we may choose a positive integer $m_n$ and an abelian
extension $L_n/K_n$ with $(\Z/p^{m_n}\Z)$-free Galois group such that
$L_n \subset L_{n+1}$ and $L = \cup L_n$.
Then one sees easily that
$$
    \Psi_{L/K}(b) = \lim_{\leftarrow}
    \Psi_{m_n,L_n/K_n}(b_n) 
$$
	for $b = (b_n) \in \mc{U}_K$.
	
Let
$$
  \mc{U}_{L/K} = \langle b = (b_n) \in \mc{U}_K \mid
  b_n \in (\mc{E}_{K_n} \cap N_{L_n/K_n} L_n^{\times})/\mc{E}_{K_n}^{p^{m_n}} \mr{\ for\ all\ } n \rangle
$$
be the group of universal norm sequences in $K$ from the extension $L/K$ (again independent
of $S$ containing the primes above $p$).  Recall that $Y_L$ denotes the
Galois group of the maximal unramified abelian pro-$p$ extension of $L$ in
which all primes above $p$ split completely.
Let $Y_{L/K}$
denote the image of the restriction map $Y_L \to Y_K$, and let
$\mc{D}_{L/K}$ be its kernel.
We remark that all primes not over $p$ must split completely in $Y_K$, so $Y_K$
is in fact the inverse limit of the $A_{K_n,S}$.  Furthermore, the image of 
$b \in \mc{U}_{L/K}$ under $\Psi_{L/K}$ lies in $I_G/I_G^2 \otimes_{\zp} Y_{L/K}$ and is independent of $S$ containing
the primes above $p$ and the ramified primes in $L/K$.  
Set $\mc{P}_{L/K} = \Psi_{L/K}(\mc{U}_{L/K})$.
We define the following quotient of $I_G Y_L/I_G^2 Y_L$:
$$
	\mc{Q}_{L/K} = I_G Y_L/I_G \mc{D}_{L/K}.
$$  
The following is verified from Theorem \ref{finclassgpiso} by taking inverse limits
(and by the above remarks, is independent of our choice of $S$).

\begin{theorem} \label{infclassgpiso}
  There is a canonical isomorphism 
  $$
    (I_G/I_G^2 \otimes_{\zp} Y_{L/K})/\mc{P}_{L/K}
    \xrightarrow{\sim} \mc{Q}_{L/K}
  $$
  of $\zp$-modules.
\end{theorem}

We end this section with the following remark.

\begin{remark}
Suppose that $L$ is Galois over some subfield $F_0$ of $F$, and set $\tgam
= \Gal(K/F_0)$.  Then the isomorphism in Theorem \ref{infclassgpiso} is an
isomorphism of $\zp[[\tgam]]$-modules.
\end{remark}

\section{Cup products with $p$} \label{results}

In this section, we combine our results on the Galois representation of Ohta with the cup product construction of the previous section.  
We return to the setting of Section \ref{galrep}.  In particular, we take $p \ge 5$, we let $F$ be the extension of $\Q(\mu_p)$ cut out by a nontrivial even character $\theta$ satisfying Hypotheses \ref{hyp} and \ref{square}, we let $K$ denote the cyclotomic $\zp$-extension of $F$, and we let $L$ be the fixed field of the kernel of $\bar{c}$ on $G_K$.  

We prove in this section that the map given by taking cup products with $p$ at the level of $K$ is surjective if and only if $U_p-1$ generates the Eisenstein ideal.
We describe the computational verification of this generation in the case that $N = 1$ and $p < 1000$.  This proves Conjecture \ref{pairingconj} for these $p$.  Not directly used in this proof but quite striking on its own, we demonstrate in Theorem \ref{identification} an identification between $U_p-1$ in the Eisenstein ideal $\I$ modulo its square and the image modulo $\tgam = \Gal(K/\Q)$-coinvariants of a universal norm in $K$ with norm $p^{-1}$ in $\Q$ under the $S$-reciprocity map $\Psi_{L/K}$ of \eqref{PsiLK}.  (Here, we have let $S_E$ consist of the primes above $p$ and any real places for any $E/\Q$.)
 
The following is a corollary of Theorem \ref{infclassgpiso}.

\begin{corollary} \label{specificcase}
	There is a canonical isomorphism of 	$\zp[[\tgam]]$-modules
  	$$
    		(I_G/I_G^2 \otimes_{\zp} 
    		Y_K)^{(1)}/\mc{P}_{L/K}^{(1)}
    		\xrightarrow{\sim} (I_G Y_L/I_G^2 Y_L)^{(1)}.
  	$$
\end{corollary}

\begin{proof}
	Note that $Y_{L/K} = Y_K$, as $L/K$ is totally ramified at $p$ by
	Lemma \ref{clemma}.  
	Furthermore, by Lemma \ref{samemodules}, we have
	$(\mc{D}_{L/K}/I_G Y_L)^{(\omega\theta^{-1})} = 0$, which implies that
	$$
		(I_G \mc{D}_{L/K}/I_G^2 Y_L)^{(1)} = 0,
	$$
	since $G \cong G^{(\theta\omega^{-1})}$.
	The result then follows from Theorem \ref{infclassgpiso}.
\end{proof}

Let $1-\zeta = (1-\zeta_{p^n}) \in \mc{U}_K$ for $\zeta = (\zeta_{p^n})$ a
generator of the Tate module of $\mu_{p^{\infty}}$.  Note that $1-\zeta$ has
norm $p^{[F:\Q(\mu_p)]}$ in $\Q$.   Let $\mf{h}$ and $\I$ be as before.
We remark that Theorem \ref{classiso} and its proof set up canonical isomorphisms
$$
	(I_G/I_G^2 \otimes_{\zp} Y_K)_{\tgam} \cong
	(I_G X_L/I_G^2 X_L)_{\tgam}
	\cong \I/\I^2,
$$
using in particular \eqref{diagram}, Lemma \ref{triveig}, and Lemma
\ref{samemodules} in obtaining the first isomorphism.

\begin{theorem} \label{identification}
    Under the canonical isomorphism 
    $$
        (I_G/I_G^2 \otimes_{\zp} Y_K)_{\tgam} \cong \I/\I^2
    $$
    the image of $\Psi_{L/K}(1-\zeta)$
    is identified with $[F:\Q(\mu_p)](1-U_p)$.
\end{theorem}

\begin{proof}
    Let $K_n = F(\mu_{p^n})$, and for each $n \ge 1$, 
    let $L_n \subset L$ be an abelian extension of $K_n$ 
    such that $\Gal(L_n/K_n)$ is $(\Z/p^{m_n}\Z)$-free for some $m_n \ge 1$, 
    with $L_n \subset L_{n+1}$ and $L = \cup L_n$. 
    Let $b = (b_n) \in \mc{U}_{L/K}$, and write 
    $$
    	b_n = N_{L_n/K_n} y_n \pmod{\mc{E}_{K_n}^{p^{m_n}}}
    $$ 
    for some $y_n \in L_n^{\times}$ for each $n \ge 1$.
    Let $J_{L_n}$ denote the ideal group of $L_n$ (which contains
    $J_{L_n,S}$ as a direct summand), and write
    $$
        y_n\mc{O}_{L_n} = \mf{B}_n \prod_{i=1}^{t_n}
        \mf{A}_{i,n}^{1-\sigma_{i,n}} \bmod
        p^{m_n}J_{L_n}
    $$
    for some $t_n \ge 0$, $\mf{A}_{i,n} \in J_{L_n,S}$ and $\sigma_{i,n}
    \in \Gal(L_n/K_n)$ with $1 \le i \le t_n$, and $\mf{B}_n$ in
    the subgroup of $J_{L_n}$ generated by the primes in
    $S_{L_n}$.  We then have (by \cite[Theorem 2.4]{mcs} and Lemma \ref{eltdef})
    $$
        \Psi_{m_n,L_n/K_n}(b_n) = \sum_{i=1}^{t_n} (\sigma_{i,n}-1)
        \otimes N_{L_n/K_n} \mf{A}_{i,n}.
    $$
    Let $\Theta_{L/K}(b)$ denote the image of $\Psi_{L/K}(b)$ in
    $(I_G X_L/I_G^2 X_L)_{\tgam}$. 
    As the image of $\Psi_{L/K}(b)$ in $(I_G Y_L/I_G^2 Y_L)_{\tgam}$
    is trivial by Corollary \ref{specificcase}, we have that
    $\Theta_{L/K}(b) = v[\mf{P}]$,
    where $[\mf{P}]$ denotes the image of any chosen
    prime $\mf{P}$ of $L$ above $p$ in
    $(I_G X_L/I_G^2 X_L)_{\tgam}$ and $v \in \zp$ denotes the sum of the
    valuations of $b$ at primes above $p$ in $K$.
    In other words, $\Theta_{L/K}(b)$ is the image in this group of the $-rv$th
    power of a geometric Frobenius $\Phi_p$ in $D_p$, where $r$
    is the residue degree of $p$ in $F/\Q$.

    Since $(\mc{P}_{L/K})_{\tgam}$ surjects onto
    the span of $U_p-1$ in $\I/\I^2$ by Theorem \ref{Sclassiso}
    and Corollary \ref{specificcase},
    there exists
    $b \in \mc{U}_{L/K}$ with $\bar{d}(\Theta_{L/K}(b))-1 = U_p-1$.
    Choosing such a $b$, we have $-rv = 1$, as $\bar{d}(\Phi_p) =
    U_p$ (again by \cite[Theorem 3.4.2]{ohta-cong}).
     Let $\epsilon$ denote the natural projection $\mc{U}_K \to \mc{U}_K^{(1)}$.
    We remark that $\mc{U}_K^{(1)}$ is free
    of rank $1$ over $\zp[[\Gamma]]$, generated by $\epsilon(1-\zeta)$.
    Since $v = -r^{-1}$ has trivial valuation at $p$, the element
    $\epsilon(b)$ also generates $\mc{U}_K^{(1)}$ 
    as a $\zp[[\Gamma]]$-module.  In particular, we have 
    $\mc{U}_K^{(1)} = \mc{U}_{L/K}^{(1)}$.  (It is also possible to see
    this using cohomological methods.)
    Since the sum of the valuations of $1-\zeta$ at primes above $p$ in $K$
    is equal to the number $g$ of such primes, we have
    $$
        \Theta_{L/K}(1-\zeta) = \Theta_{L/K}(\epsilon(1-\zeta)) =
        [\mf{P}]^{g} = \Theta_{L/K}(b)^{-rg} = \Theta_{L/K}(b)^{-[F:\Q(\mu_p)]}.
    $$
    The result follows.
\end{proof}

We require a few lemmas.
First, we employ Poitou-Tate duality to verify the following reflection principle.
Set $\mc{U}_F = \mc{E}_F \otimes \zp$. 

\begin{lemma} \label{reflection}
	We have $A_{F,S}^{(\omega^2\theta^{-1})} = 0$, and 
	$\mc{U}_F^{(\omega^2\theta^{-1})}$ is free of rank $1$ over $\mc{O}$.
\end{lemma}

\begin{proof}
	By Lemma \ref{triveigspace}, 
	we have that $A_{F,S}^{(\theta\omega^{-1})} = 0$ (or see
	\cite[Theorem 1.3.1]{lang}).  On the other hand, the 
	$(\theta\omega^{-1})$-eigenspace
	of $\mr{Br}_S(F)[p]$ is trivial by Hypothesis \ref{hyp}b. Hence, \eqref{brauer}
	implies that the $(\theta\omega^{-1})$-eigenspace of $H^2(G_{F,S},\mu_p)$ is 	trivial.  Since
	$$
		H^i(G_{F,S},\mu_p)^{(\theta\omega^{-1})} \cong
		H^i(G_{\Q,S},(\mc{O}/p\mc{O})(\omega^2\theta^{-1}))
	$$
	for any $i \ge 0$, the fact that the Euler characteristic of the latter cohomology
	groups is one implies that the $(\theta\omega^{-1})$-eigenspace of 
	$H^1(G_{F,S},\mu_p)$ is trivial as well.  By Poitou-Tate duality, this implies
	that the kernel of 
	$$
		H^2(G_{F,S},\mu_p)^{(\omega^2\theta^{-1})} \to 		
		\mr{Br}_S(F)[p]^{(\omega^2\theta^{-1})}
	$$
	is trivial.  As this kernel is the $(\omega^2\theta^{-1})$-eigenspace of 
	$A_{F,S}/p$, we have the first statement.
	
	Now, consider the $(\omega^2\theta^{-1})$-eigenspace of 
	$
		\bigoplus_{v \in S_F} H^1(G_{F_v},\mu_p).
	$
	In fact, it follows from a result of C.\ Greither (see \cite[Corollary 2.2]{me-det})
	that under Hypotheses \ref{hyp}b and c, 
	this eigenspace is a one-dimensional vector space over $\mc{O}/p\mc{O}$.	Again applying Poitou-Tate duality, the trivialities demonstrated above and 	\eqref{h1seq} imply that
	$$
		(\mc{E}_F/\mc{E}_F^p)^{(\omega^2\theta^{-1})} \cong 		H^1(G_{F,S},\mu_p)^{(\omega^2\theta^{-1})} \cong
		\bigl( \bigoplus_{v \in S_F} H^1(G_{F_v},\mu_p)
		\bigr)^{(\omega^2\theta^{-1})},
	$$
	finishing the proof.
\end{proof}

We also have the following.

\begin{lemma} \label{univnorm}
	Every element of $\mc{U}_F^{(\omega^2\theta^{-1})}$ is a 
	universal norm from $K$.
\end{lemma}

\begin{proof}
	Let $\mc{U}_{K/F}$ denote the subgroup of $\mc{U}_F$ of elements which are 	universal norms from $K$.
  	By \cite[Corollary A.2]{hs}, we have an exact sequence
  	\begin{equation} \label{sequence}
  		0 \to Y_K^{\Gamma} \otimes_{\zp} \Gamma \to \mc{U}_F/\mc{U}_{K/F}
  		\to \ker \bigl(\bigoplus_{v \in S_K} \Gamma_v \to \Gamma \bigr) 
  		\to (Y_K)_{\Gamma} \to A_{F,S} \to 0
  	\end{equation}
	(which in this case results from \cite[Theorem 6.5]{sinnott}).
  	By Lemma \ref{reflection}, we have that $A_{F,S}^{(\omega^2\theta^{-1})} = 0$.
 	Next, note that 
  	$$
  		\bigl(\bigoplus_{v \in S_K} \Gamma_v\bigr)^{(\omega^2\theta^{-1})} = 0
  	$$
  	by Hypothesis \ref{hyp}c.  Then, \eqref{sequence} implies that 
  	$(Y_K)_{\Gamma}^{(\omega^2\theta^{-1})} = 0$, 
  	so $Y_K^{(\omega^2\theta^{-1})} = 0$ by Nakayama's Lemma.  Applying 	\eqref{sequence} one more time, we see that $\mc{U}_F^{(\omega^2\theta^{-1})} 	= \mc{U}_{K/F}^{(\omega^2\theta^{-1})}$, as desired.
\end{proof}

For $q \ge 2$ even, let $\mf{h}_q(N;\mc{O})$ denote the ordinary cuspidal Hecke
algebra over $\mc{O}$ of level $Np$, weight $q$, and character $\theta\omega^{-q}$.
Let $\I_q(N;\mc{O})$ denote the Eisenstein ideal in $\mf{h}_q(N;\mc{O})$,
which is generated by $U_p-1$ and
$T_l-1-\theta(l)\kappa(l)^ql^{-1}$ for primes $l \neq p$.
Let $\mf{h}_q$ denote the localization of $\mf{h}_q(N;\mc{O})$ at the unique maximal ideal
containing $\I_q(N;\mc{O})$, and let $\I_q$ denote the resulting Eisenstein
ideal in $\mf{h}_q$.
Let
$$
    P_q = [1+p]-(1+p)^q \in \Lambda_{\mf{h}}.
$$
Then $\mf{h}_q \cong \mf{h}/P_q\mf{h}$ \cite[Corollary
3.2]{hida2}.  

\begin{lemma} \label{enough}
	We have $\I_q \cong \I/P_q\I$.
\end{lemma}

\begin{proof}
	We have
	$\I_q \cong \I/(P_q\mf{h} \,\cap\, \I)$, and we claim that
	 $P_q\mf{h}\, \cap\, \I = P_q\I$.  To see this, it suffices to
    	show that $P_q$ does not divide the characteristic ideal of the principal
    	$\Lambda_{\mf{h}}$-module $\mf{h}/\I$,
    	which we know to be generated by a power series $g_{\theta}$ with
    	$g_{\theta}((1+p)^s-1) = L_p(\theta,1-s)$ for $s \in \zp$.  Since $q$
    	is a positive even integer, the argument of Lemma \ref{tame} goes through
    	(with $\mc{O}(-1)$ replaced by $\mc{O}(1-q)$).
\end{proof}

By Lemma \ref{samemodules}, we may identify the $(\omega\theta^{-1})$-eigenspaces
of $X_K$ and $Y_K$, as well as of $A_F$ and $A_{F,S}$.  Similarly, the $(\omega\theta^{-1})$-eigenspace of $\mc{H}_K$, with $\mc{H}_K$ as in \eqref{HK}, is simply 
$Y_K^{(\omega\theta^{-1})}$
by Hypothesis \ref{hyp}b.  We are now ready to prove the following theorem on the surjectivity of pairing with $p$ under the cup product.

\begin{theorem} \label{pairgen}
  The following are equivalent:
	\begin{enumerate}
		\item[$\mr{a.}$] $\langle p, \mc{U}_K \rangle_K^{(\omega^2\theta^{-1})} 
			= X_K^{(\omega\theta^{-1})}(1)$,
		\item[$\mr{b.}$] $\langle p, \mc{E}_F/\mc{E}_F^p 			\rangle_{1,F}^{(\omega^2\theta^{-1})} 
		= A_F^{(\omega\theta^{-1})} \otimes \mu_p$,
		\item[$\mr{c.}$] $U_p-1$ generates the ideal $\I$ of $\mf{h}$,
		\item[$\mr{d.}$] $U_p-1$ generates the ideal $\I_q$ of $\mf{h}_q$.
	\end{enumerate}
\end{theorem}

\begin{proof}  
    We first consider the equivalence of conditions a and b.
   	It follows from Lemma \ref{univnorm} 
  	that the image of 	$\langle p,\mc{U}_K\rangle_K^{(\omega^2\theta^{-1})}$ under 	the natural quotient map
  	$$
    		X_K^{(\omega\theta^{-1})}(1) \to A_F^{(\omega\theta^{-1})} \otimes \mu_p
  	$$
  	equals
  	$\langle p,\mc{E}_F/\mc{E}_F^{p} \rangle_{1,F}^{(\omega^2\theta^{-1})}$ 
  	(using $p$ also to denote its image in
  	$\mc{E}_F/\mc{E}_F^{p}$).  Since this map is given by
  	taking the quotient modulo the action of the maximal ideal of
  	$\Lambda$ (again using Hypothesis \ref{hyp}c as in Lemma \ref{univnorm}), 
	Nakayama's Lemma implies that parts a and b are equivalent.

	Let $\mf{m}$ denote the maximal ideal of $\mf{h}$ containing $\I$.
	By Lemma \ref{enough} and Nakayama's Lemma, conditions c and d are both 	equivalent
  	to the image of $U_p-1$ generating $\I/\mf{m}\I$ as an $\mc{O}$-module.
  	We prove that this latter condition is equivalent to part b in a series of steps.
  	
	Let $Q$ be the quotient of $G$ modulo the $\Lambda$-submodule determined by
	the maximal ideal of $\Lambda$.
  	Let $E$ be the subextension of $L/F$ with $Q = \Gal(E/F)$.  
	It follows from the principality of $G$
  	as a $\Lambda$-module that $Q$ has dimension $1$ over $\mc{O}/p\mc{O}$.
  	By Lemma \ref{reflection}, $E$ is defined by 
  	$(\mc{E}_F/\mc{E}_F^p)^{(\omega^2\theta^{-1})}$ as 
  	a Kummer extension.  Let $\eta$ be a generator of 
  	$(\mc{E}_F/ \mc{E}_F^p)^{(\omega^2\theta^{-1})}$ as an $\mc{O}$-module.
  	Then
  	$$
  	\langle p, \mc{E}_F/\mc{E}_F^{p}\rangle_{1,F}^{(\omega^2\theta^{-1})} = 
  	\langle \eta,\mc{E}_F/\mc{E}_F^{p}\rangle_{1,F}^{(\omega^2\theta^{-1})},
  	$$
  	and both groups are necessarily generated by $\langle p,\eta \rangle_{1,F}$
  	as an $\mc{O}$-module.
  	That is, condition b is equivalent to the statement that 	
	$A_F^{(\omega\theta^{-1})} \otimes \mu_p$
  	is generated by $\langle p,\eta \rangle_{1,F}$ as an $\mc{O}$-module.
  
	Lemma \ref{triveig} and the isomorphism \eq{tensorisom} in the proof of
  	Theorem \ref{classiso} yield isomorphisms
  	\begin{equation} \label{seviso}
    		(I_Q/I_Q^2 \otimes A_F/p)^{(1)} \cong 
    		I_Q/I_Q^2 \otimes_{\mc{O}} A_F^{(\omega\theta^{-1})}/p 
    		\cong C/\mf{m}C \otimes_{\mc{O}}
    		B/\mf{m}B \cong \I/\mf{m}\I
  	\end{equation}
  	in the quotient.
 	Since $I_Q/I_Q^2$ is a one-dimensional $\mc{O}/p\mc{O}$-vector space, 
 	this sets up a noncanonical isomorphism of $\mc{O}$-modules
 	\begin{equation} \label{noncanon}
 		A_F^{(\omega\theta^{-1})} \otimes \mu_p \cong \I/\mf{m}\I.
 	\end{equation}
	If $\I/\mf{m}\I$ has $\mc{O}/p\mc{O}$-dimension greater than $1$,
	then the image of $U_p-1$ cannot
	generate $\I/\mf{m}\I$, nor can $\langle p,\eta \rangle_{1,F}$ 	generate $A_F^{(\omega\theta^{-1})} \otimes \mu_p$.
	We may therefore assume that both $\mc{O}/p\mc{O}$-vector spaces in 	\eqref{noncanon} are one-dimensional.  We are
 	reduced to proving that $\langle p,\eta \rangle_{1,F}$ is nonzero if and
 	only if the image of $U_p-1$ in $\I/\mf{m}\I$ is nonzero.
 	 
  	By Theorems \ref{infclassgpiso} and \ref{Sclassiso}, we have that
  	$$
    		(I_G/I_G^2 \otimes_{\zp} Y_K)_{\tgam}/(\mc{P}_{L/K})_{\tgam} \cong
    		(I_G Y_L/I_G^2 Y_L)_{\tgam} \cong
    		\I/(\mc{V}+\I^2),
  	$$
  	where $\mc{V}$ is the $\zp$-module generated by $U_p-1$.
  	By definition, $\mc{P}_{E/F}^{(1)}$ contains the image of
  	$(\mc{P}_{L/K})_{\tgam}$ in $(I_Q/I_Q^2 \otimes A_F/p)^{(1)}$.  
	In the proof of Theorem \ref{identification} it was shown that $\mc{U}_K^{(1)} = 	\mc{U}_{L/K}^{(1)}$, and hence the element $p$ in $F$ 
	is a universal norm from $L$.  
  	Therefore, $\mc{P}_{E/F}^{(1)}$ is in fact equal to
  	the image of $(\mc{P}_{L/K})_{\tgam}$.  
  	It follows that we have an isomorphism in the quotient
  	\begin{equation} \label{finlevel}
    		(I_Q/I_Q^2 \otimes A_F/p)^{(1)}/\mc{P}_{E/F}^{(1)}
    		\cong \I/(\mc{V}+\mf{m}\I).
  	\end{equation}
  	Now, $\langle p,\eta \rangle_{1,F}$ is nontrivial if and only if
  	$\mc{P}_{E/F}^{(1)}$ is nontrivial as well.  By \eqref{seviso} and
  	\eqref{finlevel}, the latter
  	nontriviality occurs if and only if the image of $U_p-1$ in $\I/\mf{m}\I$
  	is nonzero.
\end{proof}

We now focus on the case that $N = 1$, so $F = \Q(\mu_p)$, 
$K = \Q(\mu_{p^{\infty}})$, and $\theta = \omega^k$ for some positive
even integer $k < p$ with $p \mid B_k$ and $p \nmid B_{p+1-k}$.
In general, we call $(p,k)$ an irregular pair if $k$ is a positive even integer such that $k < p$ and $p \mid B_k$.  

The following was proven using a program
we wrote for the computational algebra package Magma.  In particular, we used its built-in routines for modular symbols (programmed by W. Stein).

\begin{theorem} \label{Upgenerates}
    For all irregular pairs $(p,k)$ with $p < 1000$,
    the element $U_p-1$ generates $\I_2$.
\end{theorem}

\noindent{\em Sketch of Proof.}  
It suffices to show that $U_p-1$ generates the abelian group $\I_2/\I_2^2$.  We did
this by computation as follows.
Let $\mc{H}_2$ denote the full cuspidal Hecke algebra over $\Z$ of weight
$2$, level $p$, and character $\omega^{k-2}$. 
We computed the
Hecke operators $U_p$ and $T_i$ in $\mc{H}_2$ for all
$i$ with $1 \le i \le (p+1)/6$
as elements of $M_d(\Z[\xi])$ for
an appropriate $d$ and a fixed primitive $(p-1)$st root of unity $\xi$, which
act on the space of cuspidal modular symbols with respect to a
fixed choice of basis. These $T_i$ generate $\mc{H}_2$ as a
$\Z$-submodule (for the case of $\Gamma_0(p)$, see \cite{as},
where the result is derived from a result of J. Sturm; the case of
$\Gamma_1(p)$ and arbitrary character follows similarly from
\cite[Corollary 9.20]{stein}).  We then considered the images
$\bar{T}_i$ and $\bar{U}_p$ of these matrices in $M_d(\Z/p^2\Z)$,
via the appropriate choice of embedding of $\Z[\xi]$ in $\zp$
followed by reduction.

Let $M$ denote the subgroup of $M_d(\Z/p^2\Z)$ spanned by the
$\bar{T}_i$.  We checked that $pM$ has $p$-rank $d$, the rank of
$\mc{H}_2$ as a $\Z$-module. We then computed the minimal integer
$N$ such that the $\bar{T}_i$ with $1 \le i \le N$ generate $M$ as
a $\Z/p^2\Z$-module.  Then $p$ and the $T_i-\sum_{0 < e \mid i}
\omega(e)^{k-2}e$ with $1 \le i \le N$ generate $\I_2$ over
$\Z_p$, provided that
$$
    (\mc{H}_2 \otimes \zp)/\I_2 \cong \mf{h}_2/\I_2 \cong \zp/B_{2,\omega^{k-2}}
$$
is of order $p$, which we checked.  We computed the images $I$ and
$J$ of $\I_2$ and $\I_2^2$, respectively, in $M$ using these
elements and their products.  We verified that $I$ is generated by
$J$ and $\bar{U}_p-1$.  Since $\I_2^2$ contains $p^2\mf{h}_2$,
this suffices to prove that $\I_2$ is generated by $\I_2^2$ and
$U_p-1$. \qed \vspace{1ex}

We remark that for all irregular pairs $(p,k)$ with $p < 1000$, we have 
$p \nmid B_{p+1-k}$. 
The following are then direct consequences of Theorems \ref{pairgen} and \ref{Upgenerates}.
  
\begin{corollary} \label{cupwithp}
    For all irregular primes $p$ with $p < 1000$, we have the
    equality
    $$
        \langle p,\mc{U}_K \rangle_K = X_K(1).
    $$
\end{corollary}

\begin{corollary} \label{cupgenerates}
    For all $p < 1000$, the pairing
    $\langle \, \cdot \, ,\, \cdot \, \rangle_{1,F}$ is surjective, and
    $\langle p,1-\zeta_p \rangle_{1,F}$
    generates its image as a $\zp[\Delta]$-module.
\end{corollary}

Thus, we have proven Theorem \ref{mainresult}
of the introduction.  Furthermore, $\langle \, \cdot \, ,\, \cdot \, \rangle_{1,F}$ is given
in each eigenspace, up to a scalar which is nonzero modulo $p$, by the formula that can now be found at {\tt www.math.mcmaster.ca/$^{\sim}$sharifi} (see
\cite[Theorem 5.1]{mcs}).  This computation has numerous interesting corollaries, in particular to the structure of unramified pro-$p$ Galois groups and to products
in the $K$-theory of cyclotomic integer rings (including $\Z$), which we intend to
discuss in later work.  For now, we refer the reader to \cite{mcs} for the
details of some of these applications 
and elaborate on just one.  

The surjectivity of the cup product for $p = 691$, which divides $B_{12}$, implies a conjecture of Ihara's on the existence of a particular relation in the $12$th graded piece of Ihara's $\Z_{691}$-Lie algebra arising from the outer Galois action on the pro-$691$ fundamental group $\pi^{[691]}$ of ${\bf P}^1-\{0,1,\infty\}$.  
We describe this conjecture briefly and refer the reader to \cite{ihara} for more details.  Specifically, if $\pi^{[p]}(i)$ denotes the $i$th term in the lower central series of $\pi^{[p]}$, one may consider the filtration of $G_{\Q}$ by the kernels $F^i G_{\Q}$ of the homomorphisms $G_{\Q} \to 
\mr{Out}(\pi^{[p]}/\pi^{[p]}(i+1))$ induced by
the outer action of $G_{\Q}$ on $\pi^{[p]}$.  The direct sum
$$
	\mf{g}_p = \bigoplus_{i = 1}^{\infty} F^i G_{\Q}/F^{i+1} G_{\Q}
$$
forms a graded $\zp$-Lie algebra with commutator induced from that of $G_{\Q}$.  
Ihara showed that the Lie algebra $\mf{g}_p$ contains nontrivial elements $\sigma_i \in \mr{gr}^i \mf{g}_p$ for odd $i \ge 3$.  P.\ Deligne conjectured that the graded $\qp$-Lie algebra $\mf{g}_p \otimes_{\zp} \qp$ is freely generated by the $\sigma_i$.  
On the other hand, Ihara conjectured a specific relation in $\mf{g}_{691}$:
$$
	2[\sigma_3,\sigma_9] - 27[\sigma_5,\sigma_7] \in 691\mr{gr}^{12} 
	\mf{g}_{691}.
$$
This is an immediate corollary of Corollary \ref{cupgenerates} for $p = 691$ and \cite[Theorem 9.11]{mcs}.  The point is that this relation is the image of a relation in the Galois group of the maximal pro-$691$ unramified outside $691$ extension of $\Q(\mu_{691})$ in which the powers of the commutators appearing in the relation are determined by values of the cup product on cyclotomic $691$-units.

\section{Selmer groups} \label{selmer}

In this section, we give one special example of the use of Theorem \ref{classiso} in the computation of Selmer groups.  That is, we consider Greenberg's Selmer group attached to Ohta's Galois module $\mf{X}$ introduced in Section \ref{prelim}.
We consider only the lattice $\mf{X}$ in the tensor product of $\mf{X}$ with the quotient field of the cuspidal Hecke algebra $\mf{h}$, and we do not introduce any cyclotomic twist.  
(Another particularly nice choice is the lattice $\mf{X}^+ \oplus \I\mf{X}^-$ twisted by the inverse of  $\det \rho$.)
  
Set $W = \mf{X} \otimes_{\mf{h}} \Hom_{\zp}(\mf{h}/\I,\qp/\zp)$.  This is isomorphic as a Hecke module to the
Pontryagin dual of $\mf{X}/\I\mf{X}$, while its Galois action is induced by that on $\mf{X}$.  Recall our fixed decomposition and inertia groups $D_p$ and $I_p$, respectively, at $p$ in $G_{\Q}$. 
We will compute the Selmer
group
$$
	\Sel(\Q,W) = \ker(H^1(G_{\Q,S},W) \to H^1(I_p,W))
$$
and the strict Selmer group
$$
	\Sel_S(\Q,W) = \ker(H^1(G_{\Q,S},W) \to H^1(D_p,W)).
$$
Such Selmer groups are defined and discussed in detail in \cite{greenberg}.  (In general, it is necessary to take a quotient of $W$ by a $D_p$-submodule in the local Selmer groups, but here that submodule is trivial as the Hodge-Tate weights specialize to $0$ and $1-k$
in weight $k$.)

We continue with the notation and hypotheses of the previous sections.
The decomposition of $\mf{X}$ as $\mf{X}^+ \oplus \mf{X}^-$
yields a decomposition of $W$ into plus and minus parts which are
$(\mf{h}/\I)$-modules.  Note that $G_{\Q}$ acts by $\phi_C$ of \eqref{phiC} on $W$.  In other words, $W$ fits in an exact sequence
\begin{equation*} \label{exactt1}
  0 \to W^+ \to W \to W/W^+ \to 0
\end{equation*}
of Galois modules, where $W^+$ is fixed by $G_{\Q}$ and
$W/W^+ = (W/W^+)^{(\omega\theta^{-1})}$.  
We have Selmer groups for $W^+$ and $W/W^+$ defined analogously to those of $W$.

We prove our result through a series of lemmas.

\begin{lemma} \label{w1+vanish}
  The Selmer group $\Sel(\Q,W^+)$ is trivial.
\end{lemma}

\begin{proof}
  We have an exact sequence
  $$
    0 \to H^1(\tgam,W^+) \to H^1(G_{\Q,S},W^+) \to
    \Hom_{\tlam}(Z_K,W^+)
  $$
  (where $Z_K$ is the Galois group of the 
  maximal abelian pro-$p$ unramified outside $p$ extension of $K$),
  and the rightmost group is trivial since $W^+$ has a trivial
  $\tgam$-action.  On the other hand, the facts that $p \nmid [F:\Q]$ and
  $K/F$ is totally ramified at primes above $p$ imply that $H^1(\tgam, W^+)$
  has trivial intersection with $\Sel(\Q, W^+)$, which must then be $0$.
\end{proof}

\begin{lemma} \label{selmerisom1}
  The natural map
  $$
    \Sel(\Q, W) \to \Sel(\Q, W/W^+)
  $$
  is an isomorphism.
\end{lemma}

\begin{proof}
  We remark that $H^2(G_{\Q,S},W^+)$ is trivial
  since $W^+$ is fixed by $\tgam$.
  As $H^2(G_{\Q,S},W^+)$ is $p$-power torsion and $G_{\Q,S}$ has
  $p$-cohomological dimension $2$, 
  this follows from the triviality of
  $$
    H^2(G_{\Q,S}, \Z/p\Z) \cong (A_{\Q(\mu_p)}/p)^{(\omega)} = 0.
  $$
  Noting also that $(W/W^+)^{(1)} = 0$, we have a commutative diagram
  $$
    \xymatrix{\UseComputerModernTips
    0 \ar[r] & H^1(G_{\Q,S}, W^+) \ar[r] \ar[d]^{\phi} & H^1(G_{\Q,S}, W)
    \ar[r] \ar[d] & H^1(G_{\Q,S}, W/W^+) \ar[d] \ar[r] & 0\\
    0 \ar[r] & H^1(I_p, W^+) \ar[r] & H^1(I_p, W) \ar[r] &
    H^1(I_p, W/W^+). &
    }
  $$
  In fact, $\phi$ factors through $H^1(D_p, W^+)$, and hence
  its image is contained in the $D_p/I_p$-invariant part of
  $H^1(I_p, W^+)$.  Let
  $\phi'$ denote the map to this group induced by $\phi$.
  By the Snake Lemma and Lemma \ref{w1+vanish}, we obtain an exact
  sequence
  $$
    0 \to \Sel(\Q, W) \to \Sel(\Q, W/W^+)
    \to \coker \phi'.
  $$
  We remark that
  $$
    H^1(I_p, W^+)^{D_p/I_p} \cong H^1(\tgam, W^+),
  $$
  which shows that $\phi'$ is surjective.
\end{proof}

\begin{lemma} \label{cohomw1}
    There is a canonical isomorphism of $\zp$-modules
    $$
        \Sel(\Q, W/W^+) \xrightarrow{\sim}
    \Hom_{\mf{h}}(\I/\I^2, W^+),
    $$
    which is given by lifting to $H^1(G_{\Q,S}, W)$
    followed by restriction to $G_{HL,S}$.  The resulting
    homomorphism has image in $W^+$ and factors through
    $(I_G X_L/I_G^2 X_L)_{\tgam} \cong \I/\I^2$, the latter isomorphism
    being as in Theorem \ref{classiso}.
\end{lemma}

\begin{proof}
  Since $(W/W^+)^{(1)} = 0$, we have an isomorphism
  $$
    H^1(G_{\Q,S},W/W^+) \cong
    H^1(G_{K,S},W/W^+)^{\tgam}.
  $$
  Hence, we have
  $$
    \Sel(\Q,W/W^+) \cong \Hom_{\tlam}(X_K,W/W^+).
  $$
  Eigenspace considerations 
  force an isomorphism
  $$
    \Hom_{\tlam}(X_K,W/W^+) \cong
    \Hom_{\Lambda}(X_K^{(\omega\theta^{-1})},W/W^+).
  $$
  The latter group is computable by Theorem \ref{blemma}:
  $$
    \Hom_{\Lambda}(X_K^{(\omega\theta^{-1})},W/W^+) \cong
    \Hom_{\mf{h}}(B/\I B,W^-).
  $$
  Since $C/\I C$ is free of rank $1$ over $\mf{h}/\I$, we have that
  \begin{equation} \label{tensoriso}
    \Hom_{\mf{h}}(B/\I B,W^-) \cong
    \Hom_{\mf{h}}(B/\I B \otimes_{\mf{h}} C/\I C, W^- \otimes_{\mf{h}} C/\I
    C).
  \end{equation}
  Since $C = \End_{\mf{h}}(\mf{X}^-,\mf{X}^+) \cong \mf{h}$, we have
  a canonical isomorphism
  $$
  	W^- \otimes_{\mf{h}} C/\I C \cong W^+.
  $$
  This and the fact that 
  $B/\I B \otimes_{\mf{h}} C/\I C \cong \I/\I^2$ now imply the existence of the desired
  isomorphism.

  Let us show that this isomorphism agrees with its description
  in the statement of this theorem.  That is, take a cocycle $f$
  with class in $\Sel(\Q,W/W^+)$, and let $\tilde{f}$ be a
  cocycle with class in $H^1(G_{\Q,S},W)$ lifting that of $f$.
  Then the restriction $\tilde{f}|_{G_{K,S}}$ lifts the homomorphism
  $f|_{G_{K,S}}$.  Furthermore, $\tilde{f}|_{G_{HL,S}}$ is a
  homomorphism with image in $W^+$.  Let $E$ denote the fixed field of the
  kernel of $\tilde{f}|_{G_{HL,S}}$, so that $\tilde{f}|_{G_{HL,S}}$
  factors through a map $g$ on $\Gal(E/HL)$.  Since
  no $\tilde{\Delta}$-invariant homomorphisms $G_{K,S} \to \Z/p\Z$ exist, 
  $\Gal(E/HL)$ is contained in
  the commutator subgroup of $\Gal(E/K)$.  As $G_{L,S}$ acts trivially
  on $W$, we have
  $$
    \tilde{f}([\sigma,\tau]) = (\sigma-1)\tilde{f}(\tau)
    = \bar{c}(\sigma)f(\tau) \in W^+.
  $$
  for $\sigma \in G_{K,S}$ and $\tau \in G_{L,S}$.  Therefore, the
  map which takes $f$ to $g$ is, by the identifications in
  Theorem \ref{classiso}, both the isomorphism in
  \eqref{tensoriso} and given by the description in the statement
  of the lemma.
\end{proof}

Lemmas \ref{selmerisom1} and \ref{cohomw1} yield the
following as an immediate corollary.

\begin{theorem} \label{selmer1}
    There is a canonical isomorphism of $\zp$-modules
    $$
        \Sel(\Q,W) \xrightarrow{\sim}
    \Hom_{\mf{h}}(\I/\I^2,W^+).
    $$
\end{theorem}

As for the strict Selmer group, we obtain the following. Recall
that $\mc{V}$ denotes the $\zp$-module generated by $U_p-1$.

\begin{theorem} \label{strictselmer1}
  There is a canonical isomorphism of $\zp$-modules
  $$
    \Sel_S(\Q,W) \xrightarrow{\sim}
    \Hom_{\mf{h}}(\I/(\mc{V}+\I^2),W^+).
  $$
\end{theorem}

\begin{proof}
  Let $\vartheta$ denote the homomorphism
  $$
    \vartheta \colon H^1(G_{\Q,S},W^+) \to H^1(D_p,W^+).
  $$
  Again, $\Sel_S(\Q,W^+) = 0$, and as in the proof of
  Theorem \ref{selmer1}, we obtain an exact sequence
  $$
    0 \to \Sel_S(\Q,W) \to \Sel_S(\Q,W/W^+) \xrightarrow{\delta}
    \coker \vartheta.
  $$
  Note that $\Sel_S(\Q,W/W^+) \cong \Sel(\Q,W/W^+)$, as 
  $Y_K^{(\omega\theta^{-1})} \cong X_K^{(\omega\theta^{-1})}$.  
  Since $D_p$ acts trivially on    
  $W^+$, we have
  $$
        H^1(D_p,W^+) \cong H^1(\tgam,W^+)
        \oplus H^1(D_p/I_p,W^+),
  $$
  and the cokernel of $\vartheta$ is canonically isomorphic to
  $H^1(D_p/I_p,W^+)$.

  By Lemma \ref{cohomw1} and the Snake Lemma, $\delta$ induces a map
  $$
    \Hom_{\mf{h}}(\I/\I^2,W^+)
    \to H^1(D_p/I_p,W^+),
  $$
  which is given by identifying
  $\I/\I^2$ and $(I_G X_L/I_G^2 X_L)_{\tgam}$ as
  in Theorem \ref{classiso},
  restricting to the image of $D_p \cap G_{L,S}$ in this group, which
  is isomorphic to a quotient of $D_p/I_p$, and finally inflating.
  By Theorem
  \ref{Sclassiso}, the kernel of this map is canonically isomorphic to
  $\Hom_{\mf{h}}(\I/(\mc{V}+\I^2),W^+)$, as desired.
\end{proof}

Let $\mf{W} = \mf{X} \otimes_{\mf{h}} \Hom_{\zp}(\mf{h},\qp/\zp)$.  We obtain as a corollary the following small piece of a ``main conjecture'' for an Eisenstein component of the Hida representation, which should be a direct analogue of \cite[Conjecture 2.2]{greenberg}.  Note that $U_p-1$ occurs as a factor in the two-variable $p$-adic $L$-function 
of Mazur and K.\ Kitagawa at the trivial character (see \cite[Theorem 4.8]{kitagawa}, for instance).

\begin{proposition}
  If $\I/\I^2$ is generated as a $\zp$-module by $U_p-1$, then
  $\Sel_S(\Q,\mf{W}) = 0$.
\end{proposition}

\begin{proof}
  Since $\I$ is a principal ideal of $\mf{h}$ generated by $U_p-1$, 
  we have an exact sequence
  $$
    0 \to W \to \mf{W} \xrightarrow{\cdot(U_p-1)} \mf{W} \to 0.
  $$
  Since the image of $c|_{I_p}$ generates $\mf{h}$
  and $d|_{I_p}$ is trivial, we have 
  $\mf{W}^{G_{\Q,S}} \subseteq \mf{W}^+$ and
  $\mf{W}^{D_p} \subseteq \mf{W}^+$.  Since $d$ has image $U_p$
  on the inverse of a Frobenius at $p$ in $D_p$ and $U_p-1$ generates $\I$,
  we have furthermore that $\mf{W}^{G_{\Q,S}}$ and $\mf{W}^{D_p}$ are
  contained in $W^+$.  Finally, since the images of $b$ and $d-1$ are contained 
  in $\I$, we have equality, i.e.,
  $$
  	\mf{W}^{G_{\Q,S}} = \mf{W}^{D_p} = W^+.
  $$
  A diagram chase then provides an exact sequence of Selmer groups
  $$
    \Sel_S(\Q,W) \to \Sel_S(\Q,\mf{W}) \xrightarrow{\cdot(U_p-1)}
    \Sel_S(\Q,\mf{W}).
  $$
  By Theorem \ref{strictselmer1}, $\Sel_S(\Q,W) = 0$, and the result
  follows from Nakayama's Lemma applied to the finitely generated 
  $\mf{h}$-module
  $\Hom_{\zp}(\Sel_S(\Q,\mf{W}),\qp/\zp)$.
\end{proof}

\appendix
\section{Eigenspaces} \label{eigenspaces}

In this appendix, we prove a decomposition result (Proposition \ref{tensor}) for the structure of eigenspaces of tensor products of $\zp[\tdel]$-modules.  We use only its consequence (Lemma \ref{triveig}) for trivial eigenspaces in this article, but the general
result would be useful should one wish to consider more general eigenspaces.

Let $\tdel$ be a finite abelian group of order prime to $p$.
Let $\tdel^*$ denote the group of $\overline{\qp}$-valued characters on $\tdel$. 
Given a character $\psi \in \tdel^*$, we let $R_{\psi}$ denote the unramified extension of $\zp$ generated by the values of $\psi$.  We denote the set of $\Gal(\overline{\qp}/\qp)$-conjugacy classes of characters in $\tdel^*$ by $\Sigma$, and we denote the conjugacy class of $\psi$ by $[\psi]$.  We have a ring decomposition (see \cite[Section 1.3]{mw})
$$
	\zp[\tdel] \cong \prod_{[\psi] \in \Sigma} R_{\psi}.
$$
In particular, the surjection $\zp[\tdel] \to R_{\psi}$ induced by $\psi$ is split.

If $A$ is $\zp[\tdel]$-module, then
we let 
$$
    A^{(\psi)} = A  \otimes_{\zp[\tdel]} R_{\psi} 
    \cong \{ a \in A \otimes_{\zp} R_{\psi} \mid \delta a = \psi(\delta)
   	a \mr{\ for\ all\ } \delta \in \tdel \},
$$
where the map $\zp[\tdel] \to R_{\psi}$ in the tensor product is the surjection
induced by $\psi$.  We remark that $A^{(\psi)}$ inherits the structure of a
$R_{\psi}[\tdel]$-module.
Again, we have a decomposition
\begin{equation} \label{decomp}
	A \cong \bigoplus_{[\psi] \in \Sigma} A^{(\psi)}.
\end{equation}
In particular $A^{(\psi)}$ and $A^{(\chi)}$ are equal as subgroups and isomorphic as $\zp[\tdel]$-submodules of $A$ if $[\psi] = [\chi]$.

Let $\mc{R}$ denote the subring of $\overline{\qp}$ generated by all character values of $\tdel$ over $\zp$.  Let $A^{\mc{R}} = A \otimes_{\zp} \mc{R}$, an $\mc{R}[\tdel]$-module.
We set 
$$
	A^{\psi} = 
	A^{\mc{R}} \otimes_{\mc{R}[\tdel]} \mc{R} \cong 
	A \otimes_{\zp[\tdel]} \mc{R} \cong A^{(\psi)} \otimes_{R_{\psi}} \mc{R},
$$
where $\zp[\tdel] \to \mc{R}$ and $\mc{R}[\tdel] \to \mc{R}$ are the maps induced by $\psi$ and $R_{\psi} \to \mc{R}$ is the natural inclusion. 
After extending scalars to $\mc{R}$, our decomposition takes the simpler form
$$
	A^{\mc{R}} \cong \bigoplus_{\chi \in \tdel^*} A^{\chi}.
$$

We wish to work with eigenspaces of tensor products.  Let us begin by proving
such a result after extension of scalars.

\begin{lemma} \label{Otensor}
	Let $\psi \in \tdel^*$, and let $A$ and $B$ be $\zp[\tdel]$-modules.  Then 
	$$
		(A \otimes_{\zp} B)^{\psi} \cong \bigoplus_{\chi \in \tdel^*}
		(A^{\chi} \otimes_{\mc{R}} B^{\psi\chi^{-1}}).
	$$
\end{lemma}

\begin{proof}
	Note that
	$$
		A \otimes_{\zp} B \otimes_{\zp} \mc{R} \cong
		A^{\mc{R}} \otimes_{\mc{R}} B^{\mc{R}}.
	$$
	Decomposing $A^{\mc{R}}$ and $B^{\mc{R}}$, we have that
	$$
		A^{\mc{R}} \otimes_{\mc{R}} B^{\mc{R}}
		\cong (\bigoplus_{\chi \in \tdel^*} A^{\chi}) \otimes_{\mc{R}}
		(\bigoplus_{\varphi \in \tdel^*} B^{\varphi})
		\cong \bigoplus_{\chi,\varphi \in \tdel^*} (A^{\chi} \otimes_{\mc{R}}
		B^{\varphi\chi^{-1}}).
	$$
	
\end{proof}

We now prove our decomposition result for eigenspaces of tensor products.

\begin{proposition} \label{tensor}
	Let $A$ and $B$ be $\zp[\tdel]$-modules, and let $\varphi \in \tdel^*$.  
	Let 
	$$
		\Sigma^2_{\varphi} = \{ ([\chi],[\psi]) \in \Sigma \times \Sigma \mid \exists \,		\theta \in [\chi] \mr{\ with\ } \varphi\theta^{-1} \in [\psi] \}.
	$$
	Then
	we have a decomposition of $\zp[\tdel]$-modules
	$$
		(A \otimes_{\zp} B)^{(\varphi)} \cong 
		\bigoplus_{([\chi],[\psi]) \in \Sigma^2_{\varphi}}
		(A^{(\chi)} \otimes_{\zp} B^{(\psi)})^{(\varphi)}.
	$$
\end{proposition}

\begin{proof}
	Since $A$ and $B$ both decompose into eigenspaces as in \eqref{decomp},
	it suffices to show that 
	$$
		(A^{(\chi)} \otimes_{\zp} B^{(\psi)})^{(\varphi)} = 0
	$$
	for $\chi, \psi \in \tdel^*$ with $([\chi],[\psi]) \notin \Sigma^2_{\varphi}$.
	Note that
	$(A^{(\chi)})^{\theta} = 0$ for $\theta \in \tdel^*$ unless $\theta \in [\chi]$,
	in which case  $(A^{(\chi)})^{\theta} \cong A^{\theta}$.	Therefore, Lemma \ref{Otensor} implies that
	$$
		(A^{(\chi)} \otimes_{\zp} B^{(\psi)})^{\varphi} \cong
		\bigoplus_{\theta \in [\chi]}
		A^{\theta} \otimes_{\mc{R}} (B^{(\psi)})^{\varphi\theta^{-1}},
	$$
	which is trivial since $\varphi\theta^{-1} \notin [\psi]$ for all $\theta \in [\chi]$.
\end{proof}

We have the following result for the special case of the trivial eigenspace of the tensor
product.

\begin{lemma} \label{triveig}
	Let $A$ and $B$ be $\zp[\tdel]$-modules.  Then
	$$
		(A^{(\chi)} \otimes_{\zp} B)^{(1)} 
		\cong (A^{(\chi)} \otimes_{\zp} B^{(\chi^{-1})})^{(1)} 
		\cong A^{(\chi)} \otimes_{R_{\chi}} B^{(\chi^{-1})}.
	$$
\end{lemma}

\begin{proof}
	Note that $\Sigma_1^2 = \{([\chi],[\chi^{-1}]) \mid [\chi] \in \Sigma\}$.
	The result then follows from Proposition \ref{tensor} and the isomorphisms
	$$
		(A^{(\chi)} \otimes_{\zp} B^{(\chi^{-1})})^{(1)}
		\cong \frac{A^{(\chi)} \otimes_{\zp} B^{(\chi^{-1})}}{\langle 
		\chi(\delta) a
		\otimes b - a \otimes \chi(\delta)b \mid a,b,\delta\rangle}
		\cong A^{(\chi)} \otimes_{R_{\chi}} B^{(\chi^{-1})},
	$$
	(where $a$, $b$, and $\delta$ in the quotient run over elements of 
	$A^{(\chi)}$, $B^{(\chi^{-1})}$ and $\tdel$, respectively).
\end{proof}

\renewcommand{\baselinestretch}{1}

\vspace{2ex} \footnotesize \noindent
Romyar Sharifi \\
Department of Mathematics and Statistics\\
McMaster University\\
1280 Main Street West\\
Hamilton, Ontario L8S 4K1, Canada\\
e-mail address: {\tt sharifi@math.mcmaster.ca}\\
web page: {\tt http://www.math.mcmaster.ca/$\thicksim$sharifi}
\end{document}